\documentclass{article}
\usepackage{colortbl,amsmath,afterpage,verbatim,url,rotating}
\usepackage[maxfloats=50]{morefloats}
\def\sp{$\spadesuit$}
\def\cl{$\clubsuit$}     
\def\di{$\diamondsuit$}
\def\he{$\heartsuit$}
\makeatletter 
\def\@biblabel#1{\hspace*{-\labelsep}}
\makeatother

\begin{document}

\title{Optimal conditional expectation \\ at the video poker game Jacks or Better}
\author{S. N. Ethier,\thanks{Department of Mathematics, University of Utah, 155 South 1400 East, Salt Lake City, UT 84112, USA. \texttt{ethier@math.utah.edu}.  Partially supported by a grant from the Simons Foundation (209632).} \ 
John Jungtae Kim,\thanks{1801-5161 Yonge Street, Toronto, Ontario, M2N 0E9, Canada. \texttt{johnnakim@gmail.com}.} \ 
and Jiyeon Lee\phantom{,}\thanks{Department of Statistics, Yeungnam University, 214-1 Daedong, Kyeongsan, Kyeongbuk 712-749, South Korea. \texttt{leejy@yu.ac.kr}.}}
\date{}
\maketitle

\begin{abstract}
There are $134{,}459$ distinct initial hands at the video poker game Jacks or Better, taking suit exchangeability into account.  A computer program can determine the optimal strategy (i.e., which cards to hold) for each such hand, but a complete list of these strategies would require a book-length manuscript.  Instead, a hand-rank table, which fits on a single page and reproduces the optimal strategy perfectly, was found for Jacks or Better as early as the mid 1990s.  Is there a systematic way to derive such a hand-rank table?  We show that there is indeed, and it involves finding the exact optimal conditional expected return, given the initial hand.  In the case of Jacks or Better (paying $800,50,25,9,6,4,3,2,1,0$), this is a random variable with $1{,}153$ distinct values, of which 766 correspond to garbage hands for which it is optimal to draw five new cards.  We describe the hands corresponding to each of the remaining 387 values of the optimal conditional expected return (sorted from largest to smallest) and show how this leads readily to an optimal strategy hand-rank table for Jacks or Better.  Of course, the method applies to other video poker games as~well.
\end{abstract}

\section{Introduction}

The video poker game Jacks or Better is a one-player game played at a video monitor.   The player places a bet.  He then receives five cards face up on the screen, with each of the $\binom{52}{5}=2{,}598{,}960$ possible hands equally likely (the order of the five cards is irrelevant).  For each card, the player must decide whether to hold that card or not.  Thus, there are $2^5=32$ ways to play the hand.  If he holds $k$ cards ($0\le k\le5$), he is dealt $5-k$ new cards, with each of the $\binom{47}{5-k}$ possibilities equally likely.  The player then receives his payout, which depends on the amount he bet and the rank of his final hand.

The pay table for full-pay Jacks or Better is shown in Table~\ref{payoff-odds}.  We use the term ``full-pay'' to emphasize the fact that there are similar versions of the game with less favorable pay tables.  (Often the payouts for full house and flush, 9 and 6, are reduced, for example to 8 and 5.  A full-pay game is also called a 9/6 game.)  It should also be mentioned that, typically, to qualify for the 800 for 1 payout on a royal flush, the player must make a maximum bet (five times the minimum bet), which we assume he does, and we regard that maximum bet amount as one unit.

\begin{table}[htb]
\caption{\label{payoff-odds}The full-pay Jacks or Better pay table, assuming a maximum bet.\medskip}
\catcode`@=\active \def@{\hphantom{0}}
\catcode`#=\active \def#{\hphantom{{,}}}
\renewcommand{\arraystretch}{1.}
\begin{center}
\begin{tabular}{lc}\hline
\noalign{\smallskip}
rank & payoff odds  \\
\noalign{\smallskip}\hline
\noalign{\smallskip}
royal flush & 800 for 1  \\
straight flush & @50 for 1  \\
four of a kind & @25 for 1  \\
full house & @@9 for 1  \\
flush & @@6 for 1  \\
straight & @@4 for 1  \\
three of a kind & @@3 for 1  \\
two pairs & @@2 for 1  \\
one pair, jacks or better & @@1 for 1  \\
other & @@0 for 1  \\
\noalign{\smallskip}\hline
\end{tabular}
\end{center}
\end{table}

The video poker game Jacks or Better dates back only to about 1979, and the first analyses were based on computer simulation or approximate computation (Frome 1989, Weber and Scruggs 1992, Gordon 1996, 2000).  Some authors (Wong 1988, Paymar 1992) argued that a slightly suboptimal strategy, being easier to implement, is preferable.  Eventually, the exactly optimal strategy was published by Dancer (1996) in the form of a one-page hand-rank table.   Dancer and Daily (2003, Table 6.2) later found a simpler, but still optimal, hand-rank table.  Ethier (2010, Section 17.1) provided a textbook treatment, using the hand-rank table of Dancer and Daily without confirming it; subsequent analysis by Kim (2012) (and perhaps others) confirmed its accuracy.  Marshall (2006) published a 357-page manual listing the correct play for every hand of the video poker game Deuces Wild; no such book has been published for Jacks or Better.

To determine the optimal strategy at Jacks or Better, it suffices, for each of the player's $2{,}598{,}960$ possible initial hands, to determine which of the 32 drawing strategies maximizes his conditional expected return.  We can reduce the amount of work required by nearly a factor of 20 by taking equivalence of initial hands into account.  Let us call two initial hands \textit{equivalent} if they have the same five denominations and if the corresponding denominations have the same suits after a permutation of (\cl,\di,\he,\sp).  For example, the equivalence class containing A\cl-A\di-A\he-K\cl-Q\di\ has 24 members because the A-K suit can be chosen in four ways, then the A-Q suit can be chosen in three ways, and finally the remaining A suit can be chosen in two ways.

How many equivalence classes are there associated with a particular set of denominations?  Let us consider the case of a hand with five distinct denominations $m_1,m_2,m_3,m_4,m_5$ with $2\le m_1<m_2<m_3<m_4<m_5\le 14$.  (Here 11--14 signify J--A.  There are $\binom{13}{5}=1{,}287$ ways to choose the denominations.)  We number the suits of denominations $m_1,m_2,m_3,m_4,m_5$ by $n_1,n_2,n_3,n_4,n_5\in\{1,2,3,4\}$.  Since we are concerned only with equivalence classes, we choose $n_1,n_2,n_3,n_4,n_5$ successively, using the smallest available integer for each suit that does not appear in a lower denomination.  Thus, $n_1=1$, $n_2\le n_1+1$, $n_3\le\max(n_1,n_2)+1$, and so on.  It is easy to see that there is a one-to-one correspondence between the set of such $(n_1,n_2,n_3,n_4,n_5)$ and the set of equivalence classes of hands with denominations $(m_1,m_2,m_3,m_4,m_5)$.  By direct enumeration (rather than by combinatorial analysis) we find that there are 51 equivalence classes.  See Table~\ref{equiv-classes} for a list, which also includes the other types of hands.

Table~\ref{equiv-classes} shows that there are exactly
\begin{equation*}
51\binom{13}{5}+20\binom{13}{1,3,9}+8\binom{13}{2,1,10}+5\binom{13}{1,2,10}+3\binom{13}{1,1,11}=134{,}459\quad
\end{equation*}
equivalence classes.  For a group-theoretic derivation of this number, see Alspach (2007).  As a check, we compute the total number of hands by summing the sizes of the equivalence classes:
\begin{eqnarray*}
&&(1\cdot4+15\cdot12+35\cdot24)\binom{13}{5}+(8\cdot12+12\cdot24)\binom{13}{1,3,9}\\
&&\qquad{}+(4\cdot12+4\cdot24)\binom{13}{2,1,10}+(1\cdot4+3\cdot12+1\cdot24)\binom{13}{1,2,10}\\
&&\qquad{}+(1\cdot4+2\cdot12)\binom{13}{1,1,11}=2{,}598{,}960. 
\end{eqnarray*}

Our program (written in C++ in 2007) methodically cycles through each of the $134{,}459$ equivalence classes.  For each one it computes the 32 conditional expected returns and determines which is largest and if it is uniquely the largest.  It stores this information in a file as it proceeds.  The program consists of a main program and two subroutines:

\begin{verbatim}
int main()
void optimal(int m0[ ], int n0[ ], int no, int sz, double t0[ ])
int poker(int l0[ ])
\end{verbatim}

\begin{table}[htb]
\caption{\label{equiv-classes}List of equivalence classes of
initial player hands in Jacks or Better, together with the size of each equivalence class. (From Ethier 2010, p.~551.)\medskip}
\catcode`@=\active \def@{\hphantom{0}}
\catcode`#=\active \def#{\hphantom{{,}}}
\tabcolsep=1.8mm
\begin{small}
\begin{center}
\begin{tabular}{ccccccccccc}\hline
\noalign{\smallskip}
\multicolumn{11}{c}{five distinct denominations $(x,y,z,v,w)$: $\binom{13}{5}=1{,}287$ ways} \\
\multicolumn{11}{c}{(includes hands ranked no pair, straight, flush, straight flush, royal flush)} \\
\noalign{\smallskip}
$(1,1,1,1,1)$ & @4 && $(1,1,2,3,3)$ & 24 && $(1,2,2,1,2)$ & 12 && $(1,2,3,2,1)$ & 24 \\
$(1,1,1,1,2)$ & 12 && $(1,1,2,3,4)$ & 24 && $(1,2,2,1,3)$ & 24 && $(1,2,3,2,2)$ & 24 \\
$(1,1,1,2,1)$ & 12 && $(1,2,1,1,1)$ & 12 && $(1,2,2,2,1)$ & 12 && $(1,2,3,2,3)$ & 24 \\
$(1,1,1,2,2)$ & 12 && $(1,2,1,1,2)$ & 12 && $(1,2,2,2,2)$ & 12 && $(1,2,3,2,4)$ & 24 \\
$(1,1,1,2,3)$ & 24 && $(1,2,1,1,3)$ & 24 && $(1,2,2,2,3)$ & 24 && $(1,2,3,3,1)$ & 24 \\
$(1,1,2,1,1)$ & 12 && $(1,2,1,2,1)$ & 12 && $(1,2,2,3,1)$ & 24 && $(1,2,3,3,2)$ & 24 \\
$(1,1,2,1,2)$ & 12 && $(1,2,1,2,2)$ & 12 && $(1,2,2,3,2)$ & 24 && $(1,2,3,3,3)$ & 24 \\
$(1,1,2,1,3)$ & 24 && $(1,2,1,2,3)$ & 24 && $(1,2,2,3,3)$ & 24 && $(1,2,3,3,4)$ & 24 \\
$(1,1,2,2,1)$ & 12 && $(1,2,1,3,1)$ & 24 && $(1,2,2,3,4)$ & 24 && $(1,2,3,4,1)$ & 24 \\
$(1,1,2,2,2)$ & 12 && $(1,2,1,3,2)$ & 24 && $(1,2,3,1,1)$ & 24 && $(1,2,3,4,2)$ & 24 \\
$(1,1,2,2,3)$ & 24 && $(1,2,1,3,3)$ & 24 && $(1,2,3,1,2)$ & 24 && $(1,2,3,4,3)$ & 24 \\
$(1,1,2,3,1)$ & 24 && $(1,2,1,3,4)$ & 24 && $(1,2,3,1,3)$ & 24 && $(1,2,3,4,4)$ & 24 \\
$(1,1,2,3,2)$ & 24 && $(1,2,2,1,1)$ & 12 && $(1,2,3,1,4)$ & 24 \\
\noalign{\smallskip}\hline
\noalign{\smallskip}
\multicolumn{11}{c}{one pair $(x,x,y,z,v)$: $\binom{13}{1,3,9}=2{,}860$ ways} \\
\noalign{\smallskip}
$(1,2,1,1,1)$ & 12 && $(1,2,1,2,3)$ & 24 && $(1,2,3,1,1)$ & 24 && $(1,2,3,3,3)$ & 12 \\
$(1,2,1,1,2)$ & 12 && $(1,2,1,3,1)$ & 24 && $(1,2,3,1,2)$ & 24 && $(1,2,3,3,4)$ & 12 \\
$(1,2,1,1,3)$ & 24 && $(1,2,1,3,2)$ & 24 && $(1,2,3,1,3)$ & 24 && $(1,2,3,4,1)$ & 24 \\
$(1,2,1,2,1)$ & 12 && $(1,2,1,3,3)$ & 24 && $(1,2,3,1,4)$ & 24 && $(1,2,3,4,3)$ & 12 \\
$(1,2,1,2,2)$ & 12 && $(1,2,1,3,4)$ & 24 && $(1,2,3,3,1)$ & 24 && $(1,2,3,4,4)$ & 12 \\
\noalign{\smallskip}\hline
\noalign{\smallskip}
\multicolumn{11}{c}{two pairs $(x,x,y,y,z)$: $\binom{13}{2,1,10}=858$ ways} \\
\noalign{\smallskip}
$(1,2,1,2,1)$ & 12 && $(1,2,1,3,1)$ & 24 && $(1,2,1,3,3)$ & 24 && $(1,2,3,4,1)$ & 12 \\
$(1,2,1,2,3)$ & 12 && $(1,2,1,3,2)$ & 24 && $(1,2,1,3,4)$ & 24 && $(1,2,3,4,3)$ & 12 \\
\noalign{\smallskip}\hline
\noalign{\smallskip}
\multicolumn{11}{c}{three of a kind $(x,x,x,y,z)$: $\binom{13}{1,2,10}=858$ ways} \\
\noalign{\smallskip}
$(1,2,3,1,1)$ & 12 && $(1,2,3,1,4)$ & 12 && $(1,2,3,4,4)$ & @4 && & \\
$(1,2,3,1,2)$ & 24 && $(1,2,3,4,1)$ & 12 && & && & \\
\noalign{\smallskip}\hline
\noalign{\smallskip}
\multicolumn{11}{c}{full house $(x,x,x,y,y)$: $\binom{13}{1,1,11}=156$ ways} \\
\noalign{\smallskip}
$(1,2,3,1,2)$ & 12 && $(1,2,3,1,4)$ & 12 && & && & \\
\noalign{\smallskip}\hline
\noalign{\smallskip}
\multicolumn{11}{c}{four of a kind $(x,x,x,x,y)$: $\binom{13}{1,1,11}=156$ ways} \\
\noalign{\smallskip}
$(1,2,3,4,1)$ & @4 && & && & && & \\
\noalign{\smallskip}\hline
\end{tabular}
\end{center}
\end{small}
\end{table}
\afterpage{\clearpage}

The subroutine \texttt{poker} returns the rank (0--9 corresponding to the 10 ranks in Table~\ref{payoff-odds}) of a hand given by a vector of five distinct integers, each in the range 1--52.  The subroutine \texttt{optimal} prints a line of output that describes the optimal strategy and conditional expected return, given a hand's five denominations, five suits, equivalence class number, and equivalence class size.  (The vector \texttt{t0} keeps a running count of the payout distribution under the optimal strategy.)  Each call of the subroutine \texttt{optimal} calls the subroutine \texttt{poker}
$$
\sum_{k=0}^5\binom{5}{k}\binom{47}{5-k}=\binom{52}{5}
$$
times.  Finally, the program \texttt{main} loops through all $1{,}287$ vectors of distinct initial denominations, and for each one calls the subroutine \texttt{optimal} 51 times (see Table~\ref{equiv-classes}); then the $2{,}860$ pairs are treated, and so on until all equivalence classes have been analyzed.  

We remark that no attempt was made to optimize the efficiency of this program, whose runtime is 3--12 hours, depending on the computer used.  Shackleford (2016) described two shortcuts that resulted in a claimed runtime of three seconds.  The first is one we have already used, suit exchangeability.  The second involves precomputing all draw probabilities.  This is a clever idea, which would be useful if one were attempting to analyze multiple pay tables.

\begin{table}[ht]
\caption{\label{output}Abbreviated output of C++ program.\medskip}
\catcode`@=\active \def@{\hphantom{0}}
\catcode`#=\active \def#{\hphantom{{,}}}
\renewcommand{\arraystretch}{1.}
\tabcolsep=.125cm
\begin{small}
\begin{center}
\begin{tabular}{l}\hline
\noalign{\smallskip}
\texttt{1,4,2,3,4,5,6,1,1,1,1,1,1,1,1,1,1,383484750,u}\\
\texttt{2,12,2,3,4,5,6,1,1,1,1,2,1,1,1,1,1,30678780,u}\\
\texttt{3,12,2,3,4,5,6,1,1,1,2,1,1,1,1,1,1,30678780,u}\\
\texttt{4,12,2,3,4,5,6,1,1,1,2,2,1,1,1,1,1,30678780,u}\\
\texttt{5,24,2,3,4,5,6,1,1,1,2,3,1,1,1,1,1,30678780,u}\\
$\vdots$\\
\texttt{134455,4,14,14,14,14,9,1,2,3,4,1,1,1,1,1,1,191742375,n}\\
\texttt{134456,4,14,14,14,14,10,1,2,3,4,1,1,1,1,1,1,191742375,n}\\
\texttt{134457,4,14,14,14,14,11,1,2,3,4,1,1,1,1,1,1,191742375,n}\\
\texttt{134458,4,14,14,14,14,12,1,2,3,4,1,1,1,1,1,1,191742375,n}\\
\texttt{134459,4,14,14,14,14,13,1,2,3,4,1,1,1,1,1,1,191742375,n}\\
\noalign{\smallskip}\hline
\end{tabular}
\end{center}
\end{small}
\end{table}

A small sample of our output is shown in Table~\ref{output}.  Notice that the output is formatted so that it can be imported into an Excel file via a \texttt{.csv} (comma-separated values) file.  The Excel file has $134{,}459$ rows and 19 columns.  The interpretation of columns 1--19 (or A--S) is as follows:\medskip

(1)  equivalence class number\par
(2)  size of equivalence class (necessarily 4, 12, or 24)\par
(3--7)  denominations of cards 1--5 in initial hand (11--14 signify J--A)\par
(8--12)  suits of cards 1--5 in initial hand (numbered 1--4; see text)\par
(13--17)  optimal strategy for cards 1--5 (1 if card is held, 0 if not)\par
(18)  expected return under optimal strategy${}\times{}7{,}669{,}695$\par
(19)  u if optimal strategy is uniquely optimal, n if not\medskip\par

The number $7{,}669{,}695$ is a common denominator (not necessarily the least one) of the various conditional probabilities, in the sense that
$$
\text{l.c.m.}\bigg\{\binom{47}{1},\binom{47}{2},\binom{47}{3},\binom{47}{4},\binom{47}{5}\bigg\}=5\binom{47}{5}=7{,}669{,}695.
$$
As was pointed out by Ethier (2010, p.~552), there are only two cases of nonuniqueness: hands with four of a kind (the odd card may be discarded or not) and hands of the form TTJQK with no more than two cards of the same suit (either ten may be discarded). Thus, it is reasonable to speak of \textit{the} optimal strategy:  it is essentially unique.  Also included in the output (but not shown in Table~\ref{output}) is the distribution of return under optimal play, with which we will not be concerned here.  See Ethier (2010, Table 17.4), for the details.

In the next section we show how to use the output in Table~\ref{output} to obtain the optimal conditional expected return at Jacks or Better, given the initial hand.

\section{Optimal conditional expectation}

The Excel file obtained from the \texttt{.csv} file mentioned above, with $134{,}459$ rows and 19 columns A--S, contains the optimal conditional expected return, but it is not in a useful form.  To make it more useful, we begin by sorting the rows, first according to column R (expected return${}\times{}7{,}669{,}695$) in descending order and then according to column A (equivalence class number) in ascending order.  We then renumber column A from 1 to $134{,}459$ and delete column S (uniqueness indicator).  Instead of Table~\ref{output} we get Table~\ref{sorted-output}.

\begin{table}[htb]
\caption{\label{sorted-output}Abbreviated list of equivalence classes of
initial player hands in Jacks or Better, sorted by optimal expected return $E$.\medskip}
\catcode`@=\active \def@{\hphantom{0}}
\catcode`#=\active \def#{\hphantom{{,}}}
\renewcommand{\arraystretch}{1.}
\tabcolsep=.125cm
\begin{small}
\begin{center}
\begin{tabular}{rrccccccccccccccccccr}\hline
\noalign{\smallskip}
no. & size && \multicolumn{5}{c}{denominations} && \multicolumn{5}{c}{suits} && \multicolumn{5}{c}{opt. strategy} & $E\times7669695$\\
\noalign{\smallskip}
\hline
\noalign{\smallskip}
1&4&&10&11&12&13&14&&1&1&1&1&1&&1&1&1&1&1&6135756000\\
2&4&&2&3&4&5&6&&1&1&1&1&1&&1&1&1&1&1&383484750\\
3&4&&2&3&4&5&14&&1&1&1&1&1&&1&1&1&1&1&383484750\\
4&4&&3&4&5&6&7&&1&1&1&1&1&&1&1&1&1&1&383484750\\
5&4&&4&5&6&7&8&&1&1&1&1&1&&1&1&1&1&1&383484750\\
$\vdots$#\\
134455&24&&5&6&7&9&10&&1&2&3&4&2&&0&0&0&0&0&2741280\\
134456&24&&5&6&8&9&10&&1&2&3&1&4&&0&0&0&0&0&2741280\\
134457&24&&5&6&8&9&10&&1&2&3&4&2&&0&0&0&0&0&2741280\\
134458&24&&5&6&7&9&10&&1&2&3&4&1&&0&0&0&0&0&2741080\\
134459&24&&5&6&8&9&10&&1&2&3&4&1&&0&0&0&0&0&2741080\\
\noalign{\smallskip}\hline
\noalign{\smallskip}
\end{tabular}
\end{center}
\end{small}
\end{table}

Call this Sheet 1.  Then, in a new worksheet in the same Excel file, called Sheet 2, we create a condensed file, with one row per distinct conditional expected return.  To do this, we first copy the contents of Sheet 1 to Sheet 2.   In Sheet 2 we insert a new column S that has a 1 or a 0 in each cell, depending on whether $E\times 7{,}669{,}695$ is smaller than in the previous row or not.  With a ``copy'' and ``paste special'' operation, we can then replace the cell formulas by just numbers.  Next we sort the rows, first according to column S (just described) in descending order and then according to column A (sorted equivalence class number) in ascending order.  We can then delete all rows except the first $1{,}153$, containing the distinct values of $E\times 7{,}669{,}695$.  (There are other ways to delete duplicates in Excel, but this works even in other spreadsheet applications such as Libre Office.)  Next we delete the 16 columns (B--Q) corresponding to equivalence class size, denominations, suits, and optimal strategy, since this information is already in Sheet 1, and we delete the last column, which now has all 1s.  This leaves only two columns, column A (smallest equivalence class number corresponding to the value of $E\times 7{,}669{,}695$ in that row) and column B ($E\times 7{,}669{,}695$). 

Now we add information to Sheet 2.  We insert a new column B with $=\text{A2}-\text{A1}$ in cell B1 and extended to all of column B.  This is the number of equivalence classes with the given value of $E$.  Then we insert another new column B, starting with $=\text{C2}+\text{B1}$ in cell B2 (1 in cell B1).  These are cumulative equivalence class numbers, useful for finding in Sheet 1 the hands corresponding to an entry in Sheet 2.  After another ``copy'' and ``paste special'' operation, we replace column A by the numbers 1 through $1{,}153$.  Next, we want to add equivalence class size information to Sheet 2.  First, we add three new columns to Sheet 1, columns S, T, and U, starting with $=\text{IF(B1}$=$4,1,0)$ in cell S1, $=\text{IF(B1}$=$12,1,0)$ in cell T1, and $=\text{IF(B1}$=$24,1,0)$ in cell U1.  Next, in new columns D, E, and F of Sheet 2, start with $=\text{SUMIF(\$Sheet1.R\$1:\$Sheet1.R\$134459}$$,\text{G1},$ $\text{\$Sheet1.S\$1:\$Sheet1.S\$134459)}$ in cell D1, and similarly in cells E1 and F1 but with S replaced by T and U.  The number of equivalence classes with given conditional expected return is in column C, and columns D, E, and F break this number down according to equivalence class size, 4, 12, or 24.  Next we insert a new column G with $=4*\text{D1}+12*\text{E1}+24*\text{F1}$ in cell G1.  This is the probability, multiplied by $\binom{52}{5}$, of the specified conditional expected return $E$.  Finally, we add a column I with a decimal approximation of $E$.  

The resulting Sheet 2 is reproduced in full as Table~\ref{387list1} (the first 387 values of $E$) and Table~\ref{766list1} (the last 766 values of $E$) in the Appendix.  Actually, we omit the cumulative equivalence class number (in Tables~\ref{387list1} and \ref{766list1}, not in the spreadsheet), which is no longer needed.  But the table still lacks one important feature.  We would like to describe the hands corresponding to each value of the conditional expected return, or at least the first 387 values.  The last 766 values correspond to garbage hands for which it is optimal to draw five new cards.  This is the only step in our algorithm that is not automated.  It therefore requires the most time (and is most prone to error).  For these descriptions, reproduced in Table~\ref{387list1}, we adopt the following notation:

\begin{quote}
Abbreviations:  RF $=$ royal flush, SF $=$ straight flush, 4K $=$ four of a kind, FH $=$ full house, F $=$ flush, S $=$ straight, 3K $=$ three of a kind, 2P $=$ two pairs, HP $=$ high pair (jacks or better), LP $=$ low pair.  $n$-RF, $n$-SF, $n$-F, and $n$-S refer to $n$-card hands that have the potential to become RF, SF, F, and S, respectively.  5-4K, 5-FH, 3-3K, 4-2P, 2-HP, and 2-LP have a slightly different meaning:  For example, 3-3K is a 3-card three of a kind, i.e., the potential is already realized (though improvement is possible).  T, J, Q, K, A denote ten, jack, queen, king, ace.  H denotes any high card (J--A).  Cards within braces are of the same suit; cards outside the braces---or within different braces---are of different suits.  $s$ is the number of straights, disregarding suits, that can be made from the hand held, and $h$ denotes the number of high cards in the hand held.  fp, sp, and 9sp denote flush penalty, straight penalty, and 9 straight penalty (explained later).   
\end{quote}

Most authors distinguish straights by the number of gaps or ``insides''.  We prefer using $s$ because it handles hands such as A234, for which $s=1$, more naturally.  Although it has no gaps, one must regard it as a one-gap hand.

We can use the Excel file, especially Sheet 2, to confirm certain known results, such as the overall expected return under optimal play.  In a new column J we multiply columns G and H.  The sum of column J ($1{,}153$ terms) is $19{,}842{,}315{,}923{,}796$, which when divided by $2{,}598{,}960\times7{,}669{,}695$ gives the desired result,
$$
\frac{1{,}653{,}526{,}326{,}983}{1{,}661{,}102{,}543{,}100}\approx0.995439043695.
$$
This is consistent with Ethier (2010, p.~552).  Thus, using the property that the expectation of a conditional expectation is the unconditional expectation, we confirm the well-known fact that Jacks or Better has expected return of 99.5439\,\% with optimal play. 

Similarly, we can readily verify that the median conditional expected return under optimal play is $4{,}452/5{,}405\approx0.823682$.  This corresponds to item 84 (2-LP) in Table~\ref{387list1}.  We might also observe that the probability of a garbage hand (one for which it is optimal to draw five new cards) is $703/21{,}658\approx0.0324591$, less than one chance in 30.

As a way of confirming the descriptions in Table~\ref{387list1} we check that the number of equivalence classes can be derived by combinatorics.   Let us illustrate this with two examples.  
\begin{enumerate}
\item Conditional expectation 80 is 4-F: $h=2$, with six categories of hands listed (separated by semi-colons).  The numbers of equivalence classes are $6\binom{8}{3}3-6=1{,}002$; $6\binom{8}{2}-1=167$; $4\binom{8}{2}3-2=334$; $3\cdot8\cdot2=48$; $3\cdot8=24$; and $6\binom{8}{2}2-2=334$.  The sum is $1{,}909$.
\item Conditional expectation 83 is 4-S: (2--8,T)TJQK (not 4-RF, 5-F, 4-SF, 3-RF, or 4-F), of which there are $7(51-14)+(20-8)=271$ equivalence classes, where we use Table~\ref{equiv-classes} to see that 14 of the 51 equivalence classes for hands of the form $(x,y,z,v,w)$, and 8 of the 20 hands of the form $(x,x,y,z,v)$, are ruled out as 4-RF, 5-F, 4-SF, 3-RF, or 4-F.
\end{enumerate}

\noindent We have confirmed all 387 cases.  

Is there a simpler way to describe optimal play?  We address this question in the next section.

\section{Optimal strategy hand-rank table}

The goal here is to derive an optimal strategy hand-rank table systematically.  This is a ranked list of descriptions of hands held, and these descriptions are not mutually exclusive.  Instead, the player applies the applicable description that ranks highest in the table, and this results in optimal play.  For example, the hand 8\cl-T\cl-J\cl-Q\cl-K\cl\ is correctly played by discarding the 8.  Thus, 4-RF should rank higher in the table than 5-F.

The key to the construction of a hand-rank table is Table~\ref{387list1} listing the 387 distinct values of the conditional expectation, sorted from largest to smallest, for which the optimal strategy involves holding at least one card.  Actually, only the sorted descriptions of the hand held and the initial hand are needed.  The conditional expectation values (as well as the probabilities and equivalence class information) are not needed in the construction.  

Let us introduce a convenient acronym.  We will call conditional expectation values ``CE values'' and conditional expectation numbers ``CE numbers''.  For example, the first entry in Table~\ref{387list1} has CE number 1 and CE value 800; the last entry has CE number 387 and CE value 0.428307.  Notice that CE values decrease as CE numbers increase.

\begin{table}[ht]
\caption{\label{prelim-hand-rank}Preliminary hand-rank table at Jacks or Better.  Notice that the first number in each row of the third column is the smallest conditional expectation number that was not previously accounted for.  This rule determines the ordering of the rows.  Caution:  This table does not reproduce the optimal strategy.}
\tabcolsep=1.8mm
\catcode`@=\active\def@{\phantom{0}}
\begin{footnotesize}
\begin{center}
\begin{tabular}{clc}
\hline
\noalign{\smallskip}
rank & hand held & conditional expectation numbers  \\
\noalign{\smallskip}
\hline
\noalign{\smallskip}
@1 & 5-RF &  1 \\
@2 & 5-SF & 2 \\
@3 & 5-4K & 3 \\
@4 & 4-RF & 4--14 \\
@5 & 5-FH & 15 \\
@6 & 5-F  & 16 \\
@7 & 3-3K & 17 \\
@8 & 5-S  & 18 \\
@9 & 4-SF: $s=2$ & 19--23 \\
10 & 4-2P & 24 \\
11 & 4-SF: $s=1$ & 25--31 \\
12 & 2-HP & 32 \\
13 & 3-RF & 33--79 \\
14 & 4-F & 80--82 \\
15 & 4-S: $s=2$, $h=3$ & 83 \\
16 & 2-LP & 84 \\
17 & 4-S: $s=2$, $h=1$ or 2 & 85--86 \\
18 & 3-SF: $s+h=4$ & 87--98 \\
19 & 4-S: $s=2$, $h=0$ & 99 \\
20 & 3-SF: $s+h=3$ & 100--110, 112--115, 117--118, 120, 127, 130, 132 \\
21 & 2-RF: JQ & 111, 116, 119, 121--126, 129, 131, \dots\ ($^1$) \\
22 & 2-RF: JK or QK  & 128, 134, 138, 144--145, 148--149, \dots\ ($^2$) \\
23 & 4-S: $s=1$, $h=4$ & 143 \\
24 & 2-RF: JA, QA, or KA & 152, 156, 160, 162, 164--167 \\
25 & 3-SF: $s+h=2$ & 168--173, 175--180, 185, 187 \\
26 & 4-S: $s=1$, $h=3$  & 174 \\
27 & 2-RF: TJ & 181, 183, 186, 189--191, 193--194, \dots\ $(^3)$ \\
28 & 3-S: JQK & 182 \\
29 & 2-S: JQ & 184, 188, 192, 195--196, 198, 201--202, 208, 215 \\
30 & 2-RF: TQ & 199, 206, 213, 221--222, 239, 247, 265, 285, 294, 329, 352 \\
31 & 2-S: JK or QK & 204, 210, 223, 232, 234, 257 \\
32 & 1-RF: J, Q, K, or A & 212, 214, 218, 224--228, 230, 233, \dots\ $(^4)$ \\
33 & 2-RF: TK & 261, 304, 380 \\
34 & 2-S: JA, QA, or KA & 269, 310 \\
35 & 3-SF: $s+h=1$ & 386--387 \\
36 & none & 388--1153 \\
\noalign{\smallskip}
\hline
\noalign{\smallskip}
\multicolumn{3}{l}{$^1$ 111, 116, 119, 121--126, 129, 131, 133, 135--137, 139--142, 146--147, 150, 153}\\
\multicolumn{3}{l}{$^2$ 128, 134, 138, 144--145, 148--149, 151, 154--155, 157--159, 161, 163}\\
\multicolumn{3}{l}{$^3$ 181, 183, 186, 189--191, 193--194, 197, 200, 203, 205, 207, 209, 211, 216--217, 219--220,}\\
\multicolumn{3}{l}{\qquad 229, 231, 250, 289}\\
\multicolumn{3}{l}{$^4$ 212, 214, 218, 224--228, 230, 233, 235--238, 240--246, 248--249, 251--256, 258--260,}\\
\multicolumn{3}{l}{\qquad 262--264, 266--268, 270--284, 286--288, 290--293, 295--303, 305--309, 311--328, 330--351,}\\
\multicolumn{3}{l}{\qquad 353--379, 381--385}\\
\end{tabular}
\end{center}
\end{footnotesize}
\end{table}\afterpage{\clearpage}

We begin by deriving Table~\ref{prelim-hand-rank}, a preliminary version of the hand-rank table, one that does not reproduce the optimal strategy perfectly (although it is a reasonable first approximation).  We simply list the hands held in order of the smallest CE number (or largest CE value) having that optimal strategy.  But before doing so, we make a few refinements.  First, 4-SF hands are not contiguous in Table~\ref{387list1}, so we separate them into two groups, those with $s=2$ and those with $s=1$.  The same is true of 4-S hands, so we classify them according to $s$ and $h$:  CE number 83 has $(s,h)=(2,3)$; CE numbers 85 and 86 have $(s,h)=(2,2)$ or $(2,1)$; CE number 99 has $(s,h)=(2,0)$; CE number 143 has $(s,h)=(1,4)$; and CE number 174 has $(s,h)=(1,3)$.  Finally, 3-SF hands are widely separated in Table~\ref{387list1}.  The first group comprises CE numbers 87--98; for these $(s,h)=(2,2)$ or $(3,1)$, which we describe as $s+h=4$.  The next group, though not exactly contiguous, comprises CE numbers 100--110, 112--115, 117--118, 120, 127, 130, and 132; for these $(s,h)=(1,2)$, $(2,1)$, or $(3,0)$, which we describe as $s+h=3$.  The next group comprises CE numbers 168--173, 175--180, 185, and 187; for these $(s,h)=(1,1)$ or $(2,0)$, which we describe as $s+h=2$.  The final group comprises CE numbers 386 and 387, with $(s,h)=(1,0)$, or $s+h=1$.  Finally, the four types of 1-RF hands in Table~\ref{387list1} are grouped together since, if an initial hand has two high cards, it is never optimal to hold just one.  With this understanding, we have 36 categories of hands held, and Table~\ref{prelim-hand-rank} sorts them according to smallest CE number (or largest CE value).  

Now we turn to the derivation of our optimal strategy hand-rank table, Table~\ref{hand-rank}, by correcting the ``errors'' in Table~\ref{prelim-hand-rank} and simplifying it slightly.  An error occurs if a hand can be played by either of two strategies and the optimal play is ranked below the suboptimal one in the table.  For example, to show that 4-RF ranks ahead of 5-F in Table~\ref{hand-rank}, we need only check that no hand that can be played as 4-RF is listed as 5-F in Table~\ref{387list1}.  In fact, the description of 5-F (CE number 16) specifically excludes 4-RF.  

\begin{table}[ht]
\caption{\label{hand-rank}Hand-rank table for optimal strategy at Jacks or Better}
\tabcolsep=1.8mm
\catcode`@=\active\def@{\phantom{0}}
\begin{footnotesize}
\begin{center}
\begin{tabular}{clcc}
\hline
\noalign{\smallskip}
rank & hand held & cond'l expected return & cond'l exp.\ nos.  \\
\noalign{\smallskip}
\hline
\noalign{\smallskip}
@1 & 5-RF & $800$ & 1 \\
@2 & 5-SF & $50$ & 2 \\
@3 & 5-4K & $25$ & 3 \\
@4 & 4-RF & $[18.361702,19.680851]$ & 4--14 \\
@5 & 5-FH & $9$ & 15 \\
@6 & 5-F  & $6$ & 16 \\
@7 & 3-3K & $4.302498$ & 17 \\
@8 & 5-S  & $4$ & 18 \\
@9 & 4-SF & $[2.340426,3.659574]$ & 19--23, 25--31 \\
10 & 4-2P & $2.595745$ & 24 \\
11 & 2-HP & $1.536540$ & 32 \\
12 & 4-F: \{(2--9)T(J--K)A\}\{(T--K)\} & $1.276596$ & 80 \\
13 & 3-RF & $[1.286772,1.532840]$ & 33--79 \\
14 & 4-F & $[1.148936,1.276596]$ & 80--82 \\
15 & 4-S: $s=2$, $h=3$ & $0.872340$ & 83 \\
16 & 2-LP &  $0.823682$ & 84 \\
17 & 4-S: $s=2$ & $[0.680851,0.808511]$ & 85--86, 99 \\
18 & 3-SF: $s+h\ge3$ & $[0.603145,0.739130]$ & 87--98, \dots\ ($^1$) \\
19 & 4-S: \{(2--7)JQ\}KA, 9\{JQ\}KA & $0.595745$ & 143 \\
20 & 2-RF: JQ & $[0.586432,0.624545]$ & 111, 116, \dots\ ($^2$) \\
21 & 4-S: $s=1$, $h=4$ & $0.595745$ & 143 \\
22 & 2-RF: (J,Q)K or (J--K)A  & $[0.557817,0.606290]$ & 128, 134, \dots\ ($^3$) \\
23 & 3-SF: $s+h=2$, no sp & $[0.533765,0.543016]$ & 168--173 \\
24 & 4-S: 9JQK, T(JQ,JK,QK)A  & $0.531915$ & 174 \\
25 & 3-SF: $s+h=2$ & $[0.506938,0.528215]$ & 175--180, 185, 187 \\
26 & 3-S: JQK & $0.515264$ & 182 \\
27 & 2-S: JQ & $[0.488375,0.509837]$ & 184, 188, \dots\ $(^4)$ \\
28 & 2-S: JK if \{TJ\} fp & $[0.483195,0.486155]$ & 223, 234 \\
29 & 2-RF: TJ & $[0.476226,0.515325]$ & 181, 183, \dots\ $(^5)$ \\
30 & 2-S: JK or QK & $[0.479494,0.494049]$ & 204, 210, \dots\ $(^6)$ \\
31 & 2-S: QA if \{TQ\} fp & 0.474314 & 310 \\
32 & 2-RF: TQ & $[0.469565,0.497071]$ & 199, 206, \dots\ $(^7)$ \\
33 & 2-S: JA, QA, or KA & $[0.474314,0.478261]$ & 269, 310 \\
34 & 1-RF: K if \{TK\} fp and 9sp & 0.459765 & 384 \\
35 & 2-RF: TK & $[0.462165,0.478816]$ & 261, 304, 380 \\
36 & 1-RF: J, Q, K, or A & $[0.458958,0.489883]$ & 212, 214, \dots\ $(^8)$ \\
37 & 3-SF: $s+h=1$ & $[0.428307,0.443108]$ & 386--387 \\
38 & none & $[0.357391,0.363385]$ & 388--1153 \\
\noalign{\smallskip}
\hline
\noalign{\smallskip}
\multicolumn{4}{l}{$^1$ 87--98, 100--110, 112--115, 117--118, 120, 127, 130, 132}\\
\multicolumn{4}{l}{$^2$ 111, 116, 119, 121--126, 129, 131, 133, 135--137, 139--142, 146--147, 150, 153}\\
\multicolumn{4}{l}{$^3$ 128, 134, 138, 144--145, 148--149, 151--152, 154--167}\\
\multicolumn{4}{l}{$^4$ 184, 188, 192, 195--196, 198, 201--202, 208, 215}\\
\multicolumn{4}{l}{$^5$ 181, 183, 186, 189--191, 193--194, 197, 200, 203, 205, 207, 209, 211, 216--217,}\\
\multicolumn{4}{l}{\qquad 219--220, 229, 231, 250, 289}\\
\multicolumn{4}{l}{$^6$ 204, 210, 223, 232, 234, 257}\\
\multicolumn{4}{l}{$^7$ 199, 206, 213, 221--222, 239, 247, 265, 285, 294, 329, 352}\\
\multicolumn{4}{l}{$^8$ 212, 214, 218, 224--228, 230, 233, 235--238, 240--246, 248--249, 251--256, 258--260,}\\
\multicolumn{4}{l}{\qquad 262--264, 266--268, 270--284, 286--288, 290--293, 295--303, 305--309, 311--328,}\\
\multicolumn{4}{l}{\qquad 330-351, 353--379, 381--385}\\
\end{tabular}
\end{center}
\end{footnotesize}
\end{table}\afterpage{\clearpage}

We follow this procedure for every rank in the table.  The first four pat hands (5-RF, 5-SF, 5-4K, 5-FH) cannot be played advantageously using lower ranked strategies.  4-RF, as we have seen, ranks ahead of 5-F and similarly 5-S.  We can rank 5-F ahead of 5-S because no hand can be played as both 5-F and 5-S.  5-F also ranks ahead of 4-SF because there are no 5-F hands listed among those for which 4-SF is optimal.  For the same reason 5-F outranks 3-RF.  Similarly, 5-S ranks ahead of 4-SF and 3-RF, as we can verify.  We have now confirmed that the first eight ranks in Table~\ref{prelim-hand-rank} are correct, so these are part of Table~\ref{hand-rank}.  In Dancer and Daily (2003, Table 6.2) and Ethier (2010, Table 17.5) some of these ranks were combined in cases where they cannot apply to the same hand simultaneously, but, as noted by Kim (2012), it is best to keep them separate so that we can include CE values in the final hand-rank table.

Hereafter we report only the discrepancies in Table~\ref{prelim-hand-rank}, not the more frequent confirmations.  The first such observation is that 4-SF: $s=2$ and 4-SF: $s=1$ can be combined as 4-SF.  Since no hand can be played as both 4-SF and 4-2P, these two ranks can appear in either order.  We put 4-SF first because its smallest CE number is 19, while 4-2P has CE number 24.   
Next we compare 3-RF (CE numbers 33--79) and 4-F (CE numbers 80--82), and we find that there are 72 (equivalence classes of) hands with CE number 80 that can be played as 3-RF, namely 4-F: (2--9)TJQA, (2--9)TJKA, (2--9)TQKA, J,Q,K odd suit; \{T\}\{(2--9)TJA\}, \{T\}\{(2--9)TQA\}, \{T\}\{(2--9)TKA\}.  To correct this ``error'' in Table~\ref{prelim-hand-rank}, we insert a new rank just above 3-RF, namely 4-F: \{(2--9)T(J--K)A\}\{(T--K)\}, followed by 3-RF and then 4-F.  Notice that our new condensed description allows a high pair, whereas CE number 80 does not.  Again, we point out that ranks are not mutually exclusive.  Such a high pair would be played as 2-HP because 2-HP outranks 3-RF and 4-F.  

Ranks 17--20 in Table~\ref{prelim-hand-rank} can be reordered and condensed to 17 combined with 19, then 18 combined with 20  in Table~\ref{hand-rank}.  We need only check that no 3-SF hand with $s+h\ge3$ can be played as 4-S with $s=2$.  We next consider ranks 21--24 in Table~\ref{prelim-hand-rank}.  Let us compare 2-RF: JQ (23 CE numbers from 111 to 153) and 4-S: $s=1$, $h=4$ (CE number 143).  Here we have another ``error'' in Table~\ref{prelim-hand-rank}, see especially CE numbers 135, 141, and 143, which contain hands that can be played both ways.  Generally, 2-RF: JQ ranks higher than 4-S: $s=1$, $h=4$, but the hands \{(2--7)JQ\}KA and 9\{JQ\}KA can be played as the former but are included in the latter.  Thus, in Table~\ref{hand-rank}, we insert a new rank just above 2-RF: JQ, namely 4-S: \{(2--7)JQ\}KA, 9\{JQ\}KA, followed by 2-RF: JQ and then 4-S: $s=1$, $h=4$.  We follow this with the two ranks 2-RF: JK or QK, and 2-RF: JA, QA, or KA, which we combine into 2-RF: (J,Q)K or (J--K)A.  Combining them is no problem because, if we had, for example, \{JK\} and \{QA\} in the same hand, the higher rank 4-S: $s=1$, $h=4$ would apply.  We can check that no hand playable as 4-S: $s=1$, $h=4$ is listed as 2-RF: (J,Q)K or (J--K)A.  This justifies Table~\ref{hand-rank} through rank 22 (corresponding to Table~\ref{prelim-hand-rank} through rank 24).  

Ranks 25 and 26 in Table~\ref{prelim-hand-rank} come next.  These are 3-SF: $s+h=2$ (CE numbers 168--173, 175--180, 185, 187) and 4-S: $s=1$, $h=3$ (CE number 174).  We notice that CE number 174 excludes hands such as \{79J\}\{Q\}\{K\}, while including
hands such as \{89Q\}\{J\}\{K\}.  The distinction is that, when these hands are played as 3-SF, the discards reduce the probability of making a straight in the latter example but not in the former.  We say then that \{89Q\}\{J\}\{K\} has a \textit{straight penalty}, while \{79J\}\{Q\}\{K\} does not.  We notice that 3-SF: $s+h=2$ can be divided into two groups: the CE numbers less than 174 have no straight penalty, while those greater than 174 have a straight penalty.  Thus, in Table~\ref{hand-rank}, we replace these two ranks by three:  3-SF: $s+h=2$, no sp; 4-S: 9JQK, T(JQ,JK,QK)A; and 3-SF: $s+h=2$.  The description of the second category seems more explicit than 4-S: $s=1$, $h=3$.

As we have noticed, it is never optimal to hold a single card when a hand has two or more high cards.  Therefore, ranks 32, 35, and 36 of Table~\ref{prelim-hand-rank} are the last three ranks of Table~\ref{hand-rank}.  It remains to properly order ranks 27--31 and 33--34, which include 3-S: JQK, three 2-RF ranks (TJ; TQ; TK), and three 2-S ranks (JQ; JK or QK; JA, QA, or KA).  
3-S: JQK comes first, then 2-S: JQ.  Next we compare 2-RF: TJ (23 CE numbers from 181 to 289) and 2-S: JK or QK (CE numbers 204, 210, 223, 232, 234, 257).  Consider hands that can be played both ways, briefly \{TJ\}K (allowing a third card to be suited with TJ).  These are CE numbers 193, 200, 207, 211, 219, 223, and 234, which are played as 2-RF: TJ unless holding TJ incurs a flush penalty (CE numbers 223 and 234).  This implies that the correct ordering is 2-S: JK if \{TJ\} fp, 2-RF: TJ, and 2-S: JK or QK. 

Next we compare 2-RF: TQ (12 CE numbers from 199 to 352) and 2-S: JA, QA, or KA (CE numbers 269 and 310).  The only hands to which both apply are of the form \{TQ\}A (allowing a third card to be suited with TQ).  These are CE numbers 221, 239, 265, 294, and 310, which are played as 2-RF: \{TQ\} unless holding TQ incurs a flush penalty (CE number 310).  This implies that the correct ordering is 2-S: QA if \{TQ\} fp, 2-RF: TQ, and 2-S: JA, QA, KA. 

Finally, we compare 2-RF: TK (CE numbers 261, 304, 380) and 1-RF: J, Q, K, or A (146 CE numbers from 212 to 385).  The only hands to which both apply are of the form \{TK\} (allowing a third card to be suited with TK).  These are CE numbers 261, 304, 380, and 384, which are played as 2-RF: TK unless holding TK incurs a flush penalty and a 9 straight penalty (CE number 384).  This implies that the correct ordering is 1-RF: K if \{TK\} fp and 9sp, 2-RF: TK, and 1-RF: J, Q, K, or A. 

Thus, Table~\ref{hand-rank}, with 38 ranks, is obtained from Table~\ref{prelim-hand-rank}.

We notice that six CE numbers appear twice in column 4 of Table~\ref{hand-rank} (namely, 80, 143, 223, 234, 310, and 384).  By breaking down column 3 of Table~\ref{387list1} in these six cases, one can readily obtain the probabilities associated with ranks 1--38.  We leave the details to the reader.

\section{Other games}

We studied Jacks or Better because it is the oldest and best-known video poker game, but the method applies to other games.  The number of distinct values of the conditional expected return is, as Table~\ref{othergames} indicates, quite sensitive to the details of the game. 
 
\begin{table}[ht]
\caption{\label{othergames}Statistics for several video poker games.  The third column gives the number of distinct values of the optimal conditional expected return.  The fourth column gives the number of distinct values of the optimal conditional expected return requiring that at least one (non-wild) card be held.
\medskip}
\catcode`@=\active \def@{\hphantom{0}}
\catcode`#=\active \def#{\hphantom{{,}}}
\renewcommand{\arraystretch}{1.}
\tabcolsep=.125cm
\begin{small}
\begin{center}
\begin{tabular}{lrrrr}\hline
\noalign{\smallskip}
                 & optimal  &  3rd & 4th \\
game (pay table) & expected &  col. & col. \\
                 & return (\%)   &       &  \\
\noalign{\smallskip}\hline
\noalign{\smallskip}
Jacks or Better $(800,50,25,9,6,4,3,2,1,0)$ &  99.5439 & $1{,}153$ & 387 \\
Double Bonus $(800,50,160,80,50,10,7,5,3,1,1,0)$  & 100.1725 & 773 & 469 \\
Joker Wild $(800,200,100,50,20,7,5,3,2,1,1,0)$    & 100.6463 & $4{,}848$ & 521  \\
Deuces Wild $(800,200,25,15,9,5,3,2,1,0)$ & 100.7620 & $8{,}903$ & 180 \\
\noalign{\smallskip}\hline
\end{tabular}
\end{center}
\end{small}
\end{table}

We have already noticed that the optimal strategy at Jacks or Better is essentially unique (meaning that the conditional payout distribution under optimal play is unique).  The same is true of Double Bonus but not of Joker Wild or Deuces Wild.  However, which optimal strategy is used in the latter two games does not affect Table~\ref{othergames}.

The reader may be surprised to see that three of the four games listed in Table~\ref{othergames} offer expected returns greater than 100\% with optimal play.  Such games can still be found, though because of concerns about advantage players, maximum bet sizes are often limited to \$1.25.

\section{Appendix}

Here we provide the exact distribution of the optimal conditional expected return at full-pay Jacks or Better, given the initial hand.  Table~\ref{387list1} gives the 387 largest values and includes mutually exclusive hand descriptions.  These correspond to cases in which at least one card is held.  Table~\ref{766list1} gives the remaining 766 values without hand descriptions.  These correspond to cases in which five new cards are drawn.

\setcounter{table}{0}
\renewcommand{\thetable}{A\arabic{table}}
\begin{sidewaystable}[ht]
\caption{\label{387list1}The 387 largest values of the optimal conditional expected return at full pay Jacks or Better, given the initial hand, page~1.  See the companion Table A2 for the remaining 766 values.\medskip}
\tabcolsep=1.0mm
\catcode`@=\active\def@{\phantom{0}}
\begin{footnotesize}
\begin{center}
\begin{tabular}{rrcrrcl}
\hline
\noalign{\smallskip}
no. & no.~of    & numbers of       & prob.                    & opt.~condl.               & opt.~condl. & hand held in optimal strategy: description of initial hand\\
    & equiv.    & equiv.~classes   & ${}\times\binom{52}{5}$  & expected@                 & expected    & (cards within braces are of the same suit; cards outside the braces\\
    & classes   & of sizes         &                          & return@@                  & return      & @---or within different braces---are of different suits)\\
    &           & 4, 12, and 24    &                          & ${}\times5\binom{47}{5}@$  & (rounded)   & \\
\noalign{\smallskip}
\hline
\noalign{\smallskip}
1&1&$(1,0,0)$&4&6135756000&800&5-RF\\
2&9&$(9,0,0)$&36&383484750&50&5-SF\\
3&156&$(156,0,0)$&624&191742375&25&5-4K\\
4&8&$(0,8,0)$&96&150946125&19.680851&4-RF: \{(2--8,T)\}\{TJQK\}\\
\hline
\end{tabular}
\end{center}
\end{footnotesize}
\end{sidewaystable}\afterpage{\clearpage}

\setcounter{table}{0}
\renewcommand{\thetable}{A\arabic{table}}
\begin{sidewaystable}[ht]
\caption{The 387 largest values of the optimal conditional expected return at full pay Jacks or Better, given the initial hand, page~2.  See the companion Table A2 for the remaining 766 values.\medskip}
\tabcolsep=1.0mm
\catcode`@=\active\def@{\phantom{0}}
\begin{footnotesize}
\begin{center}
\begin{tabular}{rrcrrcl}
\hline
\noalign{\smallskip}
no. & no.~of    & numbers of       & prob.                    & opt.~condl.               & opt.~condl. & hand held in optimal strategy: description of initial hand\\
    & equiv.    & equiv.~classes   & ${}\times\binom{52}{5}$   & expected@                 & expected    & (cards within braces are of the same suit; cards outside the braces\\
    & classes   & of sizes         &                          & return@@                  & return      & @---or within different braces---are of different suits)\\
    &           & 4, 12, and 24    &                          & ${}\times5\binom{47}{5}@$  & (rounded)   & \\
\noalign{\smallskip}
\hline
\noalign{\smallskip}
27&40&$(0,40,0)$&480&18929460&2.468085&4-SF: \{$x$\}\{89JQ\}, \{$x$\}\{8TJQ\}, \{$x$\}\{9TJK\}, \{$x$\}\{9TQK\}, not 5-S or 2-HP\\
28&8&$(0,8,0)$&96&18766275&2.446809&4-SF: \{89JQ\}\{(J,Q)\}, \{8TJQ\}\{(J,Q)\}, \{9TJK\}\{(J,K)\}, \{9TQK\}\{(Q,K)\}\\
29&88&$(0,88,0)$&1056&18439905&2.404255&4-SF: \{$x$\}\{A234\}, \{$x$\}\{A235\}, \{$x$\}\{A245\}, \{$x$\}\{A345\}, \{$x$\}\{789J\}, \{$x$\}\{78TJ\},\\
&&&&&& \quad{} \{$x$\}\{79TJ\}, \{$x$\}\{89TQ\}, not 5-S or 2-HP\\
\hline
\end{tabular}
\end{center}
\end{footnotesize}
\end{sidewaystable}\afterpage{\clearpage}

\setcounter{table}{0}
\renewcommand{\thetable}{A\arabic{table}}
\begin{sidewaystable}[ht]
\caption{The 387 largest values of the optimal conditional expected return at full pay Jacks or Better, given the initial hand, page~3.  See the companion Table A2 for the remaining 766 values.\medskip}
\tabcolsep=1.0mm
\catcode`@=\active\def@{\phantom{0}}
\begin{footnotesize}
\begin{center}
\begin{tabular}{rrcrrcl}
\hline
\noalign{\smallskip}
no. & no.~of    & numbers of       & prob.                    & opt.~condl.               & opt.~condl. & hand held in optimal strategy: description of initial hand\\
    & equiv.    & equiv.~classes   & ${}\times\binom{52}{5}$   & expected@                 & expected    & (cards within braces are of the same suit; cards outside the braces\\
    & classes   & of sizes         &                          & return@@                  & return      & @---or within different braces---are of different suits)\\
    &           & 4, 12, and 24    &                          & ${}\times5\binom{47}{5}@$  & (rounded)   & \\
\noalign{\smallskip}
\hline
\noalign{\smallskip}
51&7&$(0,7,0)$&84&11238480&1.465310&3-RF: \{(2--8)JQK\}\{A\}\\
52&6&$(0,6,0)$&72&11231385&1.464385&3-RF: \{8\}\{(2--7)TJQ\}\\
53&6&$(0,6,0)$&72&11210100&1.461610&3-RF: \{(2--7)TJQ\}\{A\}\\
54&7&$(0,7,0)$&84&11146245&1.453284&3-RF: \{T\}\{(2--8)JQK\}\\
\hline
\end{tabular}
\end{center}
\end{footnotesize}
\end{sidewaystable}\afterpage{\clearpage}

\setcounter{table}{0}
\renewcommand{\thetable}{A\arabic{table}}
\begin{sidewaystable}[ht]
\caption{The 387 largest values of the optimal conditional expected return at full pay Jacks or Better, given the initial hand, page~4.  See the companion Table A2 for the remaining 766 values.\medskip}
\tabcolsep=1.0mm
\catcode`@=\active\def@{\phantom{0}}
\begin{footnotesize}
\begin{center}
\begin{tabular}{rrcrrcl}
\hline
\noalign{\smallskip}
no. & no.~of    & numbers of       & prob.                    & opt.~condl.               & opt.~condl. & hand held in optimal strategy: description of initial hand\\
    & equiv.    & equiv.~classes   & ${}\times\binom{52}{5}$   & expected@                 & expected    & (cards within braces are of the same suit; cards outside the braces\\
    & classes   & of sizes         &                          & return@@                  & return      & @---or within different braces---are of different suits)\\
    &           & 4, 12, and 24    &                          & ${}\times5\binom{47}{5}@$  & (rounded)   & \\
\noalign{\smallskip}
\hline
\noalign{\smallskip}
77&96&$(0,48,48)$&1728&10117470&1.319149&3-RF: (2--9)J\{TQA\}, (2--9)Q\{TJA\}, (2--9)J\{TKA\}, (2--9)K\{TJA\}, (2--9)Q\{TKA\},\\
&&&&&& \quad{}(2--9)K\{TQA\} \\
78&12&$(0,6,6)$&216&10032330&1.308048&3-RF: TJ\{TQA\}, TQ\{TJA\}, TJ\{TKA\}, TK\{TJA\}, TQ\{TKA\}, TK\{TQA\}\\
79&192&$(0,192,0)$&2304&9869145&1.286772&3-RF: \{$x$\}\{$y$TJA\}, \{$x$\}\{$y$TQA\}, \{$x$\}\{$y$TKA\}, $x,y=2,3,\ldots,9$\\
\hline
\end{tabular}
\end{center}
\end{footnotesize}
\end{sidewaystable}\afterpage{\clearpage}

\setcounter{table}{0}
\renewcommand{\thetable}{A\arabic{table}}
\begin{sidewaystable}[ht]
\caption{The 387 largest values of the optimal conditional expected return at full pay Jacks or Better, given the initial hand, page~5.  See the companion Table A2 for the remaining 766 values.\medskip}
\tabcolsep=1.0mm
\catcode`@=\active\def@{\phantom{0}}
\begin{footnotesize}
\begin{center}
\begin{tabular}{rrcrrcl}
\hline
\noalign{\smallskip}
no. & no.~of    & numbers of       & prob.                    & opt.~condl.               & opt.~condl. & hand held in optimal strategy: description of initial hand\\
    & equiv.    & equiv.~classes   & ${}\times\binom{52}{5}$   & expected@                 & expected    & (cards within braces are of the same suit; cards outside the braces\\
    & classes   & of sizes         &                          & return@@                  & return      & @---or within different braces---are of different suits)\\
    &           & 4, 12, and 24    &                          & ${}\times5\binom{47}{5}@$  & (rounded)   & \\
\noalign{\smallskip}
\hline
\noalign{\smallskip}
94&2&$(0,1,1)$&36&5512815&0.718779&3-SF: \{9JQ\}KA\\
95&12&$(0,6,6)$&216&5505720&0.717854&3-SF: (2--6)\{9TJ\}K, 7\{9TJ\}A\\
96&2&$(0,1,1)$&36&5484435&0.715079&3-SF: \{9TJ\}KA\\
97&2&$(0,1,1)$&36&5420580&0.706753&3-SF: 8\{9JQ\}K\\
\hline
\end{tabular}
\end{center}
\end{footnotesize}
\end{sidewaystable}\afterpage{\clearpage}

\setcounter{table}{0}
\renewcommand{\thetable}{A\arabic{table}}
\begin{sidewaystable}[ht]
\caption{The 387 largest values of the optimal conditional expected return at full pay Jacks or Better, given the initial hand, page~6.  See the companion Table A2 for the remaining 766 values.\medskip}
\tabcolsep=1.0mm
\catcode`@=\active\def@{\phantom{0}}
\begin{footnotesize}
\begin{center}
\begin{tabular}{rrcrrcl}
\hline
\noalign{\smallskip}
no. & no.~of    & numbers of       & prob.                    & opt.~condl.               & opt.~condl. & hand held in optimal strategy: description of initial hand\\
    & equiv.    & equiv.~classes   & ${}\times\binom{52}{5}$   & expected@                 & expected    & (cards within braces are of the same suit; cards outside the braces\\
    & classes   & of sizes         &                          & return@@                  & return      & @---or within different braces---are of different suits)\\
    &           & 4, 12, and 24    &                          & ${}\times5\binom{47}{5}@$  & (rounded)   & \\
\noalign{\smallskip}
\hline
\noalign{\smallskip}
116&74&$(0,14,60)$&1608&4759799&0.620598&2-RF (JQ): $xy$8\{JQ\}, $2\le x<y\le 7$, except \{678\}\{JQ\}\\
117&8&$(0,4,4)$&144&4746555&0.618871&3-SF: \{89J\}QK, \{89J\}QA, \{8TJ\}QA, \{9TQ\}KA\\
118&84&$(0,42,42)$&1512&4739460&0.617946&3-SF: \{345\}7(J--K), \{345\}(8--T)A, \{456\}(2,8)(J--A), \{567\}(3,9)(J--A),\\
&&&&&& \quad{} \{678\}(4,T)(J--A), (2--4)\{789\}J, 5\{789\}(Q--A), (2--5)\{89T\}Q, 6\{89T\}(K,A)\\
\hline
\end{tabular}
\end{center}
\end{footnotesize}
\end{sidewaystable}\afterpage{\clearpage}

\setcounter{table}{0}
\renewcommand{\thetable}{A\arabic{table}}
\begin{sidewaystable}[ht]
\caption{The 387 largest values of the optimal conditional expected return at full pay Jacks or Better, given the initial hand, page~7.  See the companion Table A2 for the remaining 766 values.\medskip}
\tabcolsep=1.0mm
\catcode`@=\active\def@{\phantom{0}}
\begin{footnotesize}
\begin{center}
\begin{tabular}{rrcrrcl}
\hline
\noalign{\smallskip}
no. & no.~of    & numbers of       & prob.                    & opt.~condl.               & opt.~condl. & hand held in optimal strategy: description of initial hand\\
    & equiv.    & equiv.~classes   & ${}\times\binom{52}{5}$   & expected@                 & expected    & (cards within braces are of the same suit; cards outside the braces\\
    & classes   & of sizes         &                          & return@@                  & return      & @---or within different braces---are of different suits)\\
    &           & 4, 12, and 24    &                          & ${}\times5\binom{47}{5}@$  & (rounded)   & \\
\noalign{\smallskip}
\hline
\noalign{\smallskip}
140&12&$(0,6,6)$&216&4579113&0.597040&2-RF (JQ): 89\{(2--7)JQ\}\\
141&5&$(0,1,4)$&108&4574383&0.596423&2-RF (JQ): 8\{JQ\}KA\\
142&60&$(0,30,30)$&1080&4571545&0.596053&2-RF (JQ): $x$T\{$y$JQ\}, $y$T\{$x$JQ\}, $2\le x<y\le 7$\\
143&253&$(0,37,216)$&5628&4569180&0.595745&4-S: (2--8)JQKA, not 4-RF, 5-F, 3-RF, 4-F, or 3-SF, except (2--8)\{JQ\}KA;\\
\hline
\end{tabular}
\end{center}
\end{footnotesize}
\end{sidewaystable}\afterpage{\clearpage}

\setcounter{table}{0}
\renewcommand{\thetable}{A\arabic{table}}
\begin{sidewaystable}[ht]
\caption{The 387 largest values of the optimal conditional expected return at full pay Jacks or Better, given the initial hand, page~8.  See the companion Table A2 for the remaining 766 values.\medskip}
\tabcolsep=1.0mm
\catcode`@=\active\def@{\phantom{0}}
\begin{footnotesize}
\begin{center}
\begin{tabular}{rrcrrcl}
\hline
\noalign{\smallskip}
no. & no.~of    & numbers of       & prob.                    & opt.~condl.               & opt.~condl. & hand held in optimal strategy: description of initial hand\\
    & equiv.    & equiv.~classes   & ${}\times\binom{52}{5}$   & expected@                 & expected    & (cards within braces are of the same suit; cards outside the braces\\
    & classes   & of sizes         &                          & return@@                  & return      & @---or within different braces---are of different suits)\\
    &           & 4, 12, and 24    &                          & ${}\times5\binom{47}{5}@$  & (rounded)   & \\
\noalign{\smallskip}
\hline
\noalign{\smallskip}
161&196&$(0,98,98)$&3528&4410725&0.575085&2-RF (JK, QK): $y$Q\{$x$JK\}, $x$Q\{$y$JK\}, $y$J\{$x$QK\}, $x$J\{$y$QK\}, $2\le x<y\le 8$;\\
&&&&&& \quad{} 9\{(2--8)JK\}A, 9\{(2--8)QK\}A\\
162&226&$(0,42,184)$&4920&4405995&0.574468&2-RF (JA, QA, KA): (2--9)TJ\{QA\}, (2--9)TQ\{JA\}, (2--9)TJ\{KA\}, (2--9)TK\{JA\},\\
&&&&&& \quad{} (2--9)TQ\{KA\}, (2--9)TK\{QA\}, except 9TJ\{QA\}, 9TQ\{JA\}, \{8TJ\}\{QA\},\\
\hline
\end{tabular}
\end{center}
\end{footnotesize}
\end{sidewaystable}\afterpage{\clearpage}

\setcounter{table}{0}
\renewcommand{\thetable}{A\arabic{table}}
\begin{sidewaystable}[ht]
\caption{The 387 largest values of the optimal conditional expected return at full pay Jacks or Better, given the initial hand, page~9.  See the companion Table A2 for the remaining 766 values.\medskip}
\tabcolsep=1.0mm
\catcode`@=\active\def@{\phantom{0}}
\begin{footnotesize}
\begin{center}
\begin{tabular}{rrcrrcl}
\hline
\noalign{\smallskip}
no. & no.~of    & numbers of       & prob.                    & opt.~condl.               & opt.~condl. & hand held in optimal strategy: description of initial hand\\
    & equiv.    & equiv.~classes   & ${}\times\binom{52}{5}$   & expected@                 & expected    & (cards within braces are of the same suit; cards outside the braces\\
    & classes   & of sizes         &                          & return@@                  & return      & @---or within different braces---are of different suits)\\
    &           & 4, 12, and 24    &                          & ${}\times5\binom{47}{5}@$  & (rounded)   & \\
\noalign{\smallskip}
\hline
\noalign{\smallskip}
173&63&$(0,0,63)$&1512&4093815&0.533765&3-SF: (\{234\},\{235\},\{245\})(JQ,JK,QK), (\{78T\},\{79T\})(QK,QA,KA),\\
&&&&&& \quad{} (\{346\},\{356\},\{457\},\{467\},\{568\},\{578\},\{679\},\{689\})(JQ,JK,JA,QK,QA,KA),\\
&&&&&& \quad{} J,Q,K,A distinct suits\\
174&471&$(0,0,471)$&11304&4079625&0.531915&4-S: (2--8)9JQK, (2--8)TJQA, (2--9)TJKA, (2--9)TQKA, J,Q,K,A distinct suits,\\
\hline
\end{tabular}
\end{center}
\end{footnotesize}
\end{sidewaystable}\afterpage{\clearpage}

\setcounter{table}{0}
\renewcommand{\thetable}{A\arabic{table}}
\begin{sidewaystable}[ht]
\caption{The 387 largest values of the optimal conditional expected return at full pay Jacks or Better, given the initial hand, page~10.  See the companion Table A2 for the remaining 766 values.\medskip}
\tabcolsep=1.0mm
\catcode`@=\active\def@{\phantom{0}}
\begin{footnotesize}
\begin{center}
\begin{tabular}{rrcrrcl}
\hline
\noalign{\smallskip}
no. & no.~of    & numbers of       & prob.                    & opt.~condl.               & opt.~condl. & hand held in optimal strategy: description of initial hand\\
    & equiv.    & equiv.~classes   & ${}\times\binom{52}{5}$   & expected@                 & expected    & (cards within braces are of the same suit; cards outside the braces\\
    & classes   & of sizes         &                          & return@@                  & return      & @---or within different braces---are of different suits)\\
    &           & 4, 12, and 24    &                          & ${}\times5\binom{47}{5}@$  & (rounded)   & \\
\noalign{\smallskip}
\hline
\noalign{\smallskip}
184&670&$(0,120,550)$&14640&3910291&0.509837&2-S (JQ):  $xyz$JQ, $2\le x<y<z\le 7$, J,Q distinct suits, not 4-F, except\\
&&&&&& \quad{} \{$xyz$\}\{J\}\{Q\} if $z-x\le3$\\
185&16&$(0,8,8)$&288&3909345&0.509713&3-SF:  $(x-1)$\{$xyz$\}$(z+1)$, $3\le x<y<z\le 9$, $z-x=3$\\
186&49&$(0,9,40)$&1068&3891844&0.507431&2-RF (TJ):  $xy$8\{TJ\}, $2\le x<y\le 6$, except \{568\}\{TJ\}\\
\hline
\end{tabular}
\end{center}
\end{footnotesize}
\end{sidewaystable}\afterpage{\clearpage}

\setcounter{table}{0}
\renewcommand{\thetable}{A\arabic{table}}
\begin{sidewaystable}[ht]
\caption{The 387 largest values of the optimal conditional expected return at full pay Jacks or Better, given the initial hand, page~11.  See the companion Table A2 for the remaining 766 values.\medskip}
\tabcolsep=1.0mm
\catcode`@=\active\def@{\phantom{0}}
\begin{footnotesize}
\begin{center}
\begin{tabular}{rrcrrcl}
\hline
\noalign{\smallskip}
no. & no.~of    & numbers of       & prob.                    & opt.~condl.               & opt.~condl. & hand held in optimal strategy: description of initial hand\\
    & equiv.    & equiv.~classes   & ${}\times\binom{52}{5}$   & expected@                 & expected    & (cards within braces are of the same suit; cards outside the braces\\
    & classes   & of sizes         &                          & return@@                  & return      & @---or within different braces---are of different suits)\\
    &           & 4, 12, and 24    &                          & ${}\times5\binom{47}{5}@$  & (rounded)   & \\
\noalign{\smallskip}
\hline
\noalign{\smallskip}
204&2344&$(0,415,1929)$&51276&3789203&0.494049&2-S (JK, QK):  $xyz$(J,Q)K, $2\le x<y<z\le 8$, J,Q,K distinct suits, not 4-F, except\\
&&&&&& \quad{} \{$xyz$\}\{(J,Q)\}\{K\} if $z-x\le 3$, (2--6)\{78J\}K\\
205&5&$(0,1,4)$&108&3787784&0.493864&2-RF (TJ):  78\{TJ\}A\\
206&72&$(0,12,60)$&1584&3782108&0.493124&2-RF (TQ):  $xy$8\{TQ\}, $2\le x<y\le 7$, except \{568\}\{TQ\}, \{578\}\{TQ\}, \{678\}\{TQ\}\\
\hline
\end{tabular}
\end{center}
\end{footnotesize}
\end{sidewaystable}\afterpage{\clearpage}

\setcounter{table}{0}
\renewcommand{\thetable}{A\arabic{table}}
\begin{sidewaystable}[ht]
\caption{The 387 largest values of the optimal conditional expected return at full pay Jacks or Better, given the initial hand, page~12.  See the companion Table A2 for the remaining 766 values.\medskip}
\tabcolsep=1.0mm
\catcode`@=\active\def@{\phantom{0}}
\begin{footnotesize}
\begin{center}
\begin{tabular}{rrcrrcl}
\hline
\noalign{\smallskip}
no. & no.~of    & numbers of       & prob.                    & opt.~condl.               & opt.~condl. & hand held in optimal strategy: description of initial hand\\
    & equiv.    & equiv.~classes   & ${}\times\binom{52}{5}$   & expected@                 & expected    & (cards within braces are of the same suit; cards outside the braces\\
    & classes   & of sizes         &                          & return@@                  & return      & @---or within different braces---are of different suits)\\
    &           & 4, 12, and 24    &                          & ${}\times5\binom{47}{5}@$  & (rounded)   & \\
\noalign{\smallskip}
\hline
\noalign{\smallskip}
224&109&$(0,0,109)$&2616&3726982&0.485936&1-RF (J):  $xy$78\{J\}, $2\le x<y\le 6$, $(x,y)\ne(5,6)$, not 4-SF or 4-F, except \\
&&&&&& \quad{} \{(2--4)\}\{578\}\{J\}, \{(2--4)\}\{678\}\{J\}, \{457\}\{8\}\{J\}, \{467\}\{8\}\{J\}\\
225&123&$(0,0,123)$&2952&3724230&0.485577&1-RF (J):  $xyz$9\{J\}, $2\le x<y<z\le 6$, not 4-F, except \{$xyz$\}\{9\}\{J\} if $z-x\le 3$\\
226&123&$(0,0,123)$&2952&3715974&0.484501&1-RF (J):  $xy$79\{J\}, $2\le x<y\le 6$, not 4-SF or 4-F, except \{(2--5)\}\{679\}\{J\},\\
\hline
\end{tabular}
\end{center}
\end{footnotesize}
\end{sidewaystable}\afterpage{\clearpage}

\setcounter{table}{0}
\renewcommand{\thetable}{A\arabic{table}}
\begin{sidewaystable}[ht]
\caption{The 387 largest values of the optimal conditional expected return at full pay Jacks or Better, given the initial hand, page~13.  See the companion Table A2 for the remaining 766 values.\medskip}
\tabcolsep=1.0mm
\catcode`@=\active\def@{\phantom{0}}
\begin{footnotesize}
\begin{center}
\begin{tabular}{rrcrrcl}
\hline
\noalign{\smallskip}
no. & no.~of    & numbers of       & prob.                    & opt.~condl.               & opt.~condl. & hand held in optimal strategy: description of initial hand\\
    & equiv.    & equiv.~classes   & ${}\times\binom{52}{5}$   & expected@                 & expected    & (cards within braces are of the same suit; cards outside the braces\\
    & classes   & of sizes         &                          & return@@                  & return      & @---or within different braces---are of different suits)\\
    &           & 4, 12, and 24    &                          & ${}\times5\binom{47}{5}@$  & (rounded)   & \\
\noalign{\smallskip}
\hline
\noalign{\smallskip}
240&229&$(0,0,229)$&5496&3700494&0.482483&1-RF (Q):  $xyz$8\{Q\}, $2\le x<y<z\le 7$, $(x,y,z)­\ne(5,6,7)$, not 4-SF or 4-F, except\\
&&&&&& \quad{} \{$xyz$\}\{8\}\{Q\} if $z-x\le 3$ (9 cases), \{$x$\}\{$yz$8\}\{Q\} if $y\ge5$ (9 cases)\\
241&129&$(0,0,129)$&3096&3696710&0.481989&1-RF (J):  $xy$8T\{J\}, $2\le x<y\le 6$, not 4-F, except \{568\}\{T\}\{J\}\\
242&147&$(0,27,120)$&3204&3692668&0.481462&1-RF (J):  $xy$8\{$z$J\}, $xz$8\{$y$J\}, $yz$8\{$x$J\}, $2\le x<y<z\le 6$, except \{568\}\{(2--4)J\}\\
\hline
\end{tabular}
\end{center}
\end{footnotesize}
\end{sidewaystable}\afterpage{\clearpage}

\setcounter{table}{0}
\renewcommand{\thetable}{A\arabic{table}}
\begin{sidewaystable}[ht]
\caption{The 387 largest values of the optimal conditional expected return at full pay Jacks or Better, given the initial hand, page~14.  See the companion Table A2 for the remaining 766 values.\medskip}
\tabcolsep=1.0mm
\catcode`@=\active\def@{\phantom{0}}
\begin{footnotesize}
\begin{center}
\begin{tabular}{rrcrrcl}
\hline
\noalign{\smallskip}
no. & no.~of    & numbers of       & prob.                    & opt.~condl.               & opt.~condl. & hand held in optimal strategy: description of initial hand\\
    & equiv.    & equiv.~classes   & ${}\times\binom{52}{5}$   & expected@                 & expected    & (cards within braces are of the same suit; cards outside the braces\\
    & classes   & of sizes         &                          & return@@                  & return      & @---or within different braces---are of different suits)\\
    &           & 4, 12, and 24    &                          & ${}\times5\binom{47}{5}@$  & (rounded)   & \\
\noalign{\smallskip}
\hline
\noalign{\smallskip}
261&162&$(0,22,140)$&3624&3672372&0.478816&2-RF (TK):  $xyz$\{TK\}, $2\le x<y<z\le 8$, except \{$xyz$\}\{TK\} if $z-x\le 3$\\
262&50&$(0,10,40)$&1080&3671684&0.478726&1-RF (J):  $xy$9\{7J\}, $2\le x<y\le 6$\\
263&54&$(0,27,27)$&972&3670996&0.478637&1-RF (J):  $xy$\{$z$7J\}, $xz$\{$y$7J\}, $yz$\{$x$7J\}, $2\le x<y<z\le 6$, $(x,y,z)­\ne(4,5,6)$\\
264&150&$(0,30,120)$&3240&3670652&0.478592&1-RF (J):  $xy$T\{$z$J\}, $xz$T\{$y$J\}, $yz$T\{$x$J\}, $2\le x<y<z\le 6$\\
\hline
\end{tabular}
\end{center}
\end{footnotesize}
\end{sidewaystable}\afterpage{\clearpage}

\setcounter{table}{0}
\renewcommand{\thetable}{A\arabic{table}}
\begin{sidewaystable}[ht]
\caption{The 387 largest values of the optimal conditional expected return at full pay Jacks or Better, given the initial hand, page~15.  See the companion Table A2 for the remaining 766 values.\medskip}
\tabcolsep=1.0mm
\catcode`@=\active\def@{\phantom{0}}
\begin{footnotesize}
\begin{center}
\begin{tabular}{rrcrrcl}
\hline
\noalign{\smallskip}
no. & no.~of    & numbers of       & prob.                    & opt.~condl.               & opt.~condl. & hand held in optimal strategy: description of initial hand\\
    & equiv.    & equiv.~classes   & ${}\times\binom{52}{5}$   & expected@                 & expected    & (cards within braces are of the same suit; cards outside the braces\\
    & classes   & of sizes         &                          & return@@                  & return      & @---or within different braces---are of different suits)\\
    &           & 4, 12, and 24    &                          & ${}\times5\binom{47}{5}@$  & (rounded)   & \\
\noalign{\smallskip}
\hline
\noalign{\smallskip}
279&276&$(0,48,228)$&6048&3657924&0.476932&1-RF (Q):  $xy$8\{$z$Q\}, $xz$8\{$y$Q\}, $yz$8\{$x$Q\}, $2\le x<y<z\le 7$, $(x,y,z)\ne(5,6,7)$,\\
&&&&&& \quad except (\{568\},\{578\},\{678\})\{(2--4)Q\}\\
280&20&$(0,4,16)$&432&3657236&0.476842&1-RF (J):  (2--5)89\{7J\}\\
281&86&$(0,10,76)$&1944&3656204&0.476708&1-RF (Q):  $xyz$\{8Q\}, $2\le x<y<z\le 7$, $(x,y,z)\ne(5,6,7)$, except \{$xyz$\}\{8Q\}\\
\hline
\end{tabular}
\end{center}
\end{footnotesize}
\end{sidewaystable}\afterpage{\clearpage}

\setcounter{table}{0}
\renewcommand{\thetable}{A\arabic{table}}
\begin{sidewaystable}[ht]
\caption{The 387 largest values of the optimal conditional expected return at full pay Jacks or Better, given the initial hand, page~16.  See the companion Table A2 for the remaining 766 values.\medskip}
\tabcolsep=1.0mm
\catcode`@=\active\def@{\phantom{0}}
\begin{footnotesize}
\begin{center}
\begin{tabular}{rrcrrcl}
\hline
\noalign{\smallskip}
no. & no.~of    & numbers of       & prob.                    & opt.~condl.               & opt.~condl. & hand held in optimal strategy: description of initial hand\\
    & equiv.    & equiv.~classes   & ${}\times\binom{52}{5}$   & expected@                 & expected    & (cards within braces are of the same suit; cards outside the braces\\
    & classes   & of sizes         &                          & return@@                  & return      & @---or within different braces---are of different suits)\\
    &           & 4, 12, and 24    &                          & ${}\times5\binom{47}{5}@$  & (rounded)   & \\
\noalign{\smallskip}
\hline
\noalign{\smallskip}
297&25&$(0,5,20)$&540&3646228&0.475407&1-RF (J):  (2--6)8T\{7J\}\\
298&100&$(0,20,80)$&2160&3645884&0.475362&1-RF (J):  $x$9T\{$y$J\}, $y$9T\{$x$J\}, $2\le x<y\le 6$\\
299&60&$(0,30,30)$&1080&3645540&0.475317&1-RF (J):  $xy$\{$z$9J\}, $xz$\{$y$9J\}, $yz$\{$x$9J\}, $2\le x<y<z\le 6$\\
300&25&$(0,5,20)$&540&3644508&0.475183&1-RF (J)  (2--6)7T\{8J\}\\
\hline
\end{tabular}
\end{center}
\end{footnotesize}
\end{sidewaystable}\afterpage{\clearpage}

\setcounter{table}{0}
\renewcommand{\thetable}{A\arabic{table}}
\begin{sidewaystable}[ht]
\caption{The 387 largest values of the optimal conditional expected return at full pay Jacks or Better, given the initial hand, page~17.  See the companion Table A2 for the remaining 766 values.\medskip}
\tabcolsep=1.0mm
\catcode`@=\active\def@{\phantom{0}}
\begin{footnotesize}
\begin{center}
\begin{tabular}{rrcrrcl}
\hline
\noalign{\smallskip}
no. & no.~of    & numbers of       & prob.                    & opt.~condl.               & opt.~condl. & hand held in optimal strategy: description of initial hand\\
    & equiv.    & equiv.~classes   & ${}\times\binom{52}{5}$   & expected@                 & expected    & (cards within braces are of the same suit; cards outside the braces\\
    & classes   & of sizes         &                          & return@@                  & return      & @---or within different braces---are of different suits)\\
    &           & 4, 12, and 24    &                          & ${}\times5\binom{47}{5}@$  & (rounded)   & \\
\noalign{\smallskip}
\hline
\noalign{\smallskip}
311&40&$(0,20,20)$&720&3637284&0.474241&1-RF (J):  $x$7\{$y$9J\}, $y$7\{$x$9J\}, $2\le x<y\le 6$\\
312&70&$(0,14,56)$&1512&3636940&0.474196&1-RF (Q):  $xy$9\{8Q\}, $2\le x<y\le 7$, $(x,y)­\ne(6,7)$\\
313&300&$(0,60,240)$&6480&3635908&0.474062&1-RF (Q):  $xy$T\{$z$Q\}, $xz$T\{$y$Q\}, $yz$T\{$x$Q\}, $2\le x<y<z\le 7$\\
314&556&$(0,0,556)$&13344&3635478&0.474006&1-RF (K, A):  $xyz$T\{A\}, $2\le x<y\le 5$, $6\le z\le 9$, not 4-F, except\\
\hline
\end{tabular}
\end{center}
\end{footnotesize}
\end{sidewaystable}\afterpage{\clearpage}

\setcounter{table}{0}
\renewcommand{\thetable}{A\arabic{table}}
\begin{sidewaystable}[ht]
\caption{The 387 largest values of the optimal conditional expected return at full pay Jacks or Better, given the initial hand, page~18.  See the companion Table A2 for the remaining 766 values.\medskip}
\tabcolsep=1.0mm
\catcode`@=\active\def@{\phantom{0}}
\begin{footnotesize}
\begin{center}
\begin{tabular}{rrcrrcl}
\hline
\noalign{\smallskip}
no. & no.~of    & numbers of       & prob.                    & opt.~condl.               & opt.~condl. & hand held in optimal strategy: description of initial hand\\
    & equiv.    & equiv.~classes   & ${}\times\binom{52}{5}$   & expected@                 & expected    & (cards within braces are of the same suit; cards outside the braces\\
    & classes   & of sizes         &                          & return@@                  & return      & @---or within different braces---are of different suits)\\
    &           & 4, 12, and 24    &                          & ${}\times5\binom{47}{5}@$  & (rounded)   & \\
\noalign{\smallskip}
\hline
\noalign{\smallskip}
334&195&$(0,75,120)$&3780&3615956&0.471460&1-RF (Q):  $xy$T\{9Q\}, $2\le x<y\le 7$; $x$9\{$yz$Q\}, $y$9\{$xz$Q\}, $z$9\{$xy$Q\}, $2\le x<y<z\le 7$\\
335&10&$(0,5,5)$&180&3615268&0.471371&1-RF (J):  8T\{(2--6)7J\}\\
336&20&$(0,10,10)$&360&3614924&0.471326&1-RF (J):  9T\{$xy$J\}, $2\le x<y\le 6$\\
337&10&$(0,5,5)$&180&3613548&0.471146&1-RF (J):  7T\{(2--6)8J\}\\
\hline
\end{tabular}
\end{center}
\end{footnotesize}
\end{sidewaystable}\afterpage{\clearpage}

\setcounter{table}{0}
\renewcommand{\thetable}{A\arabic{table}}
\begin{sidewaystable}[ht]
\caption{The 387 largest values of the optimal conditional expected return at full pay Jacks or Better, given the initial hand, page~19.  See the companion Table A2 for the remaining 766 values.\medskip}
\tabcolsep=1.0mm
\catcode`@=\active\def@{\phantom{0}}
\begin{footnotesize}
\begin{center}
\begin{tabular}{rrcrrcl}
\hline
\noalign{\smallskip}
no. & no.~of    & numbers of       & prob.                    & opt.~condl.               & opt.~condl. & hand held in optimal strategy: description of initial hand\\
    & equiv.    & equiv.~classes   & ${}\times\binom{52}{5}$   & expected@                 & expected    & (cards within braces are of the same suit; cards outside the braces\\
    & classes   & of sizes         &                          & return@@                  & return      & @---or within different braces---are of different suits)\\
    &           & 4, 12, and 24    &                          & ${}\times5\binom{47}{5}@$  & (rounded)   & \\
\noalign{\smallskip}
\hline
\noalign{\smallskip}
353&760&$(0,148,612)$&16464&3601164&0.469532&1-RF (K, A):  $xy$T\{$z$K\}, $xz$T\{$y$K\}, $yz$T\{$x$K\}, $2\le x<y<z\le 8$, except\\
&&&&&& \quad{} \{78T\}\{(2--6)K\}; $xy$T\{$z$A\}, $xz$T\{$y$A\}, $2\le x\le 5$, $6\le y<z\le 9$\\
354&108&$(0,12,96)$&2448&3599444&0.469307&1-RF (A):  $yz$T\{$x$A\}, $2\le x\le 5$, $6\le y<z\le 9$, except \{(78,79,89)T\}\{(2--5)A\} \phantom{JQK}\\
355&60&$(0,0,60)$&1440&3597724&0.469083&1-RF (A):  $xyz$\{$v$A\}, $2\le x<y<z\le 5$, $6\le v\le 9$, $(x,y,z,v)\ne­(3,4,5,6)$,\\
\hline
\end{tabular}
\end{center}
\end{footnotesize}
\end{sidewaystable}\afterpage{\clearpage}

\setcounter{table}{0}
\renewcommand{\thetable}{A\arabic{table}}
\begin{sidewaystable}[ht]
\caption{The 387 largest values of the optimal conditional expected return at full pay Jacks or Better, given the initial hand, page~20.  See the companion Table A2 for the remaining 766 values.\medskip}
\tabcolsep=1.0mm
\catcode`@=\active\def@{\phantom{0}}
\begin{footnotesize}
\begin{center}
\begin{tabular}{rrcrrcl}
\hline
\noalign{\smallskip}
no. & no.~of    & numbers of       & prob.                    & opt.~condl.               & opt.~condl. & hand held in optimal strategy: description of initial hand\\
    & equiv.    & equiv.~classes   & ${}\times\binom{52}{5}$   & expected@                 & expected    & (cards within braces are of the same suit; cards outside the braces\\
    & classes   & of sizes         &                          & return@@                  & return      & @---or within different braces---are of different suits)\\
    &           & 4, 12, and 24    &                          & ${}\times5\binom{47}{5}@$  & (rounded)   & \\
\noalign{\smallskip}
\hline
\noalign{\smallskip}
370&280&$(0,140,140)$&5040&3571236&0.465629&1-RF (A):  $xz$\{$yv$A\}, $xv$\{$yz$A\}, $yz$\{$xv$A\}, $yv$\{$xz$A\}, $2\le x<y\le 5$, $6\le z<v\le 9$,\\
&&&&&& \quad{} $(x,y,z,v)­\ne(4,5,6,7)$\\
371&258&$(0,129,129)$&4644&3570204&0.465495&1-RF (K, A):  $x$T\{$yz$K\}, $y$T\{$xz$K\}, $z$T\{$xy$K\}, $2\le x<y<z\le 8$;  $x$T\{$yz$A\},\\
&&&&&& \quad{} $2\le x\le 5$, $6\le y<z\le 9$\\
\hline
\end{tabular}
\end{center}
\end{footnotesize}
\end{sidewaystable}\afterpage{\clearpage}

\setcounter{table}{1}
\renewcommand{\thetable}{A\arabic{table}}
\begin{sidewaystable}[ht]
\caption{\label{766list1}The 766 smallest values of the optimal conditional expected return at full pay Jacks or Better, given the initial hand, page~1.  In each case the optimal strategy is to draw five new cards.  See Table A1 for the 387 largest values.\medskip}
\tabcolsep=1.0mm
\catcode`@=\active\def@{\phantom{0}}
\begin{footnotesize}
\begin{center}
\begin{tabular}{rrcrrccrrcrrc}
\hline
\noalign{\smallskip}
no. & no.~of    & numbers of       & prob.       & opt.~condl.          & opt.~condl.  &@@@@@@&  no. & no.~of    & numbers of       & prob.       & opt.~condl.          & opt.~condl. \\
    & equiv.    & equiv.~classes   & ${}\times\binom{52}{5}$ & expected@ & expected     &&      & equiv.    & equiv.~classes   & ${}\times\binom{52}{5}$ & expected@ & expected  \\
    & classes   & of sizes         &             & return@@             & return       &&      & classes   & of sizes         &             & return@@             & return \\
    &           & 4, 12, and 24    &         & ${}\times5\binom{47}{5}@$ & (rounded)    &&      &           & 4, 12, and 24    &         & ${}\times5\binom{47}{5}@$ & (rounded) \\
\noalign{\smallskip}
\hline
\noalign{\smallskip}
388&1&$(0,1,0)$&12&2787050&0.363385&&414&1&$(0,1,0)$&12&2778650&0.362290\\
389&2&$(0,2,0)$&24&2786650&0.363333&&415&1&$(0,0,1)$&24&2778540&0.362275\\
390&1&$(0,1,0)$&12&2783970&0.362983&&416&2&$(0,2,0)$&24&2778450&0.362263\\
391&2&$(0,2,0)$&24&2783890&0.362973&&417&2&$(0,0,2)$&48&2778420&0.362260\\
\hline
\end{tabular}
\end{center}
\end{footnotesize}
\end{sidewaystable}\afterpage{\clearpage}

\setcounter{table}{1}
\renewcommand{\thetable}{A\arabic{table}}
\begin{sidewaystable}[ht]
\caption{The 766 smallest values of the optimal conditional expected return at full pay Jacks or Better, given the initial hand, page~2.  In each case the optimal strategy is to draw five new cards.  See Table A1 for the 387 largest values.\medskip}
\tabcolsep=1.0mm
\catcode`@=\active\def@{\phantom{0}}
\begin{footnotesize}
\begin{center}
\begin{tabular}{rrcrrccrrcrrc}
\hline
\noalign{\smallskip}
no. & no.~of    & numbers of       & prob.       & opt.~condl.          & opt.~condl.  &@@@@@@&  no. & no.~of    & numbers of       & prob.       & opt.~condl.          & opt.~condl. \\
    & equiv.    & equiv.~classes   & ${}\times\binom{52}{5}$ & expected@ & expected     &&      & equiv.    & equiv.~classes   & ${}\times\binom{52}{5}$ & expected@ & expected  \\
    & classes   & of sizes         &             & return@@             & return       &&      & classes   & of sizes         &             & return@@             & return \\
    &           & 4, 12, and 24    &         & ${}\times5\binom{47}{5}@$ & (rounded)    &&      &           & 4, 12, and 24    &         & ${}\times5\binom{47}{5}@$ & (rounded) \\
\noalign{\smallskip}
\hline
\noalign{\smallskip}
440&2&$(0,0,2)$&48&2776460&0.362004&&466&3&$(0,1,2)$&60&2774740&0.361780\\
441&1&$(0,1,0)$&12&2776420&0.361999&&467&5&$(0,0,5)$&120&2774620&0.361764\\
442&1&$(0,0,1)$&24&2776300&0.361983&&468&2&$(0,2,0)$&24&2774610&0.361763\\
443&5&$(0,5,0)$&60&2776290&0.361982&&469&7&$(0,2,5)$&144&2774540&0.361754\\
\hline
\end{tabular}
\end{center}
\end{footnotesize}
\end{sidewaystable}\afterpage{\clearpage}

\setcounter{table}{1}
\renewcommand{\thetable}{A\arabic{table}}
\begin{sidewaystable}[ht]
\caption{The 766 smallest values of the optimal conditional expected return at full pay Jacks or Better, given the initial hand, page~3.  In each case the optimal strategy is to draw five new cards.  See Table A1 for the 387 largest values.\medskip}
\tabcolsep=1.0mm
\catcode`@=\active\def@{\phantom{0}}
\begin{footnotesize}
\begin{center}
\begin{tabular}{rrcrrccrrcrrc}
\hline
\noalign{\smallskip}
no. & no.~of    & numbers of       & prob.       & opt.~condl.          & opt.~condl.  &@@@@@@&  no. & no.~of    & numbers of       & prob.       & opt.~condl.          & opt.~condl. \\
    & equiv.    & equiv.~classes   & ${}\times\binom{52}{5}$ & expected@ & expected     &&      & equiv.    & equiv.~classes   & ${}\times\binom{52}{5}$ & expected@ & expected  \\
    & classes   & of sizes         &             & return@@             & return       &&      & classes   & of sizes         &             & return@@             & return \\
    &           & 4, 12, and 24    &         & ${}\times5\binom{47}{5}@$ & (rounded)    &&      &           & 4, 12, and 24    &         & ${}\times5\binom{47}{5}@$ & (rounded) \\
\noalign{\smallskip}
\hline
\noalign{\smallskip}
492&4&$(0,0,4)$&96&2772860&0.361535&&518&6&$(0,0,6)$&144&2771660&0.361378\\
493&2&$(0,2,0)$&24&2772850&0.361533&&519&2&$(0,2,0)$&24&2771650&0.361377\\
494&2&$(0,0,2)$&48&2772820&0.361529&&520&2&$(0,0,2)$&48&2771540&0.361362\\
495&1&$(0,0,1)$&24&2772750&0.361520&&521&2&$(0,0,2)$&48&2771460&0.361352\\
\hline
\end{tabular}
\end{center}
\end{footnotesize}
\end{sidewaystable}\afterpage{\clearpage}

\setcounter{table}{1}
\renewcommand{\thetable}{A\arabic{table}}
\begin{sidewaystable}[ht]
\caption{The 766 smallest values of the optimal conditional expected return at full pay Jacks or Better, given the initial hand, page~4.  In each case the optimal strategy is to draw five new cards.  See Table A1 for the 387 largest values.\medskip}
\tabcolsep=1.0mm
\catcode`@=\active\def@{\phantom{0}}
\begin{footnotesize}
\begin{center}
\begin{tabular}{rrcrrccrrcrrc}
\hline
\noalign{\smallskip}
no. & no.~of    & numbers of       & prob.       & opt.~condl.          & opt.~condl.  &@@@@@@&  no. & no.~of    & numbers of       & prob.       & opt.~condl.          & opt.~condl. \\
    & equiv.    & equiv.~classes   & ${}\times\binom{52}{5}$ & expected@ & expected     &&      & equiv.    & equiv.~classes   & ${}\times\binom{52}{5}$ & expected@ & expected  \\
    & classes   & of sizes         &             & return@@             & return       &&      & classes   & of sizes         &             & return@@             & return \\
    &           & 4, 12, and 24    &         & ${}\times5\binom{47}{5}@$ & (rounded)    &&      &           & 4, 12, and 24    &         & ${}\times5\binom{47}{5}@$ & (rounded) \\
\noalign{\smallskip}
\hline
\noalign{\smallskip}
544&1&$(0,0,1)$&24&2770460&0.361222&&570&4&$(0,4,0)$&48&2769570&0.361106\\
545&4&$(0,0,4)$&96&2770420&0.361216&&571&10&$(0,9,1)$&132&2769540&0.361102\\
546&17&$(0,0,17)$&408&2770380&0.361211&&572&7&$(0,1,6)$&156&2769500&0.361096\\
547&7&$(0,7,0)$&84&2770300&0.361201&&573&2&$(0,0,2)$&48&2769470&0.361093\\
\hline
\end{tabular}
\end{center}
\end{footnotesize}
\end{sidewaystable}\afterpage{\clearpage}

\setcounter{table}{1}
\renewcommand{\thetable}{A\arabic{table}}
\begin{sidewaystable}[ht]
\caption{The 766 smallest values of the optimal conditional expected return at full pay Jacks or Better, given the initial hand, page~5.  In each case the optimal strategy is to draw five new cards.  See Table A1 for the 387 largest values.\medskip}
\tabcolsep=1.0mm
\catcode`@=\active\def@{\phantom{0}}
\begin{footnotesize}
\begin{center}
\begin{tabular}{rrcrrccrrcrrc}
\hline
\noalign{\smallskip}
no. & no.~of    & numbers of       & prob.       & opt.~condl.          & opt.~condl.  &@@@@@@&  no. & no.~of    & numbers of       & prob.       & opt.~condl.          & opt.~condl. \\
    & equiv.    & equiv.~classes   & ${}\times\binom{52}{5}$ & expected@ & expected     &&      & equiv.    & equiv.~classes   & ${}\times\binom{52}{5}$ & expected@ & expected  \\
    & classes   & of sizes         &             & return@@             & return       &&      & classes   & of sizes         &             & return@@             & return \\
    &           & 4, 12, and 24    &         & ${}\times5\binom{47}{5}@$ & (rounded)    &&      &           & 4, 12, and 24    &         & ${}\times5\binom{47}{5}@$ & (rounded) \\
\noalign{\smallskip}
\hline
\noalign{\smallskip}
596&1&$(0,0,1)$&24&2768870&0.361014&&622&3&$(0,0,3)$&72&2767470&0.360832\\
597&2&$(0,0,2)$&48&2768860&0.361013&&623&12&$(0,6,6)$&216&2767420&0.360825\\
598&1&$(0,0,1)$&24&2768820&0.361008&&624&3&$(0,3,0)$&36&2767410&0.360824\\
599&6&$(0,0,6)$&144&2768460&0.360961&&625&1&$(0,0,1)$&24&2767380&0.360820\\
\hline
\end{tabular}
\end{center}
\end{footnotesize}
\end{sidewaystable}\afterpage{\clearpage}

\setcounter{table}{1}
\renewcommand{\thetable}{A\arabic{table}}
\begin{sidewaystable}[ht]
\caption{The 766 smallest values of the optimal conditional expected return at full pay Jacks or Better, given the initial hand, page~6.  In each case the optimal strategy is to draw five new cards.  See Table A1 for the 387 largest values.\medskip}
\tabcolsep=1.0mm
\catcode`@=\active\def@{\phantom{0}}
\begin{footnotesize}
\begin{center}
\begin{tabular}{rrcrrccrrcrrc}
\hline
\noalign{\smallskip}
no. & no.~of    & numbers of       & prob.       & opt.~condl.          & opt.~condl.  &@@@@@@&  no. & no.~of    & numbers of       & prob.       & opt.~condl.          & opt.~condl. \\
    & equiv.    & equiv.~classes   & ${}\times\binom{52}{5}$ & expected@ & expected     &&      & equiv.    & equiv.~classes   & ${}\times\binom{52}{5}$ & expected@ & expected  \\
    & classes   & of sizes         &             & return@@             & return       &&      & classes   & of sizes         &             & return@@             & return \\
    &           & 4, 12, and 24    &         & ${}\times5\binom{47}{5}@$ & (rounded)    &&      &           & 4, 12, and 24    &         & ${}\times5\binom{47}{5}@$ & (rounded) \\
\noalign{\smallskip}
\hline
\noalign{\smallskip}
648&3&$(0,0,3)$&72&2766750&0.360738&&674&7&$(0,2,5)$&144&2765700&0.360601\\
649&2&$(0,0,2)$&48&2766740&0.360737&&675&2&$(0,0,2)$&48&2765670&0.360597\\
650&1&$(0,0,1)$&24&2766710&0.360733&&676&6&$(0,3,3)$&108&2765660&0.360596\\
651&4&$(0,4,0)$&48&2766660&0.360726&&677&3&$(0,0,3)$&72&2765590&0.360587\\
\hline
\end{tabular}
\end{center}
\end{footnotesize}
\end{sidewaystable}\afterpage{\clearpage}

\setcounter{table}{1}
\renewcommand{\thetable}{A\arabic{table}}
\begin{sidewaystable}[ht]
\caption{The 766 smallest values of the optimal conditional expected return at full pay Jacks or Better, given the initial hand, page~7.  In each case the optimal strategy is to draw five new cards.  See Table A1 for the 387 largest values.\medskip}
\tabcolsep=1.0mm
\catcode`@=\active\def@{\phantom{0}}
\begin{footnotesize}
\begin{center}
\begin{tabular}{rrcrrccrrcrrc}
\hline
\noalign{\smallskip}
no. & no.~of    & numbers of       & prob.       & opt.~condl.          & opt.~condl.  &@@@@@@&  no. & no.~of    & numbers of       & prob.       & opt.~condl.          & opt.~condl. \\
    & equiv.    & equiv.~classes   & ${}\times\binom{52}{5}$ & expected@ & expected     &&      & equiv.    & equiv.~classes   & ${}\times\binom{52}{5}$ & expected@ & expected  \\
    & classes   & of sizes         &             & return@@             & return       &&      & classes   & of sizes         &             & return@@             & return \\
    &           & 4, 12, and 24    &         & ${}\times5\binom{47}{5}@$ & (rounded)    &&      &           & 4, 12, and 24    &         & ${}\times5\binom{47}{5}@$ & (rounded) \\
\noalign{\smallskip}
\hline
\noalign{\smallskip}
700&9&$(0,0,9)$&216&2764950&0.360503&&726&8&$(0,0,8)$&192&2764230&0.360409\\
701&6&$(0,0,6)$&144&2764940&0.360502&&727&9&$(0,0,9)$&216&2764220&0.360408\\
702&1&$(0,0,1)$&24&2764910&0.360498&&728&17&$(0,0,17)$&408&2764190&0.360404\\
703&5&$(0,5,0)$&60&2764900&0.360497&&729&3&$(0,0,3)$&72&2764150&0.360399\\
\hline
\end{tabular}
\end{center}
\end{footnotesize}
\end{sidewaystable}\afterpage{\clearpage}

\setcounter{table}{1}
\renewcommand{\thetable}{A\arabic{table}}
\begin{sidewaystable}[ht]
\caption{The 766 smallest values of the optimal conditional expected return at full pay Jacks or Better, given the initial hand, page~8.  In each case the optimal strategy is to draw five new cards.  See Table A1 for the 387 largest values.\medskip}
\tabcolsep=1.0mm
\catcode`@=\active\def@{\phantom{0}}
\begin{footnotesize}
\begin{center}
\begin{tabular}{rrcrrccrrcrrc}
\hline
\noalign{\smallskip}
no. & no.~of    & numbers of       & prob.       & opt.~condl.          & opt.~condl.  &@@@@@@&  no. & no.~of    & numbers of       & prob.       & opt.~condl.          & opt.~condl. \\
    & equiv.    & equiv.~classes   & ${}\times\binom{52}{5}$ & expected@ & expected     &&      & equiv.    & equiv.~classes   & ${}\times\binom{52}{5}$ & expected@ & expected  \\
    & classes   & of sizes         &             & return@@             & return       &&      & classes   & of sizes         &             & return@@             & return \\
    &           & 4, 12, and 24    &         & ${}\times5\binom{47}{5}@$ & (rounded)    &&      &           & 4, 12, and 24    &         & ${}\times5\binom{47}{5}@$ & (rounded) \\
\noalign{\smallskip}
\hline
\noalign{\smallskip}
752&14&$(0,14,0)$&168&2763380&0.360299&&778&4&$(0,4,0)$&48&2762180&0.360142\\
753&8&$(0,0,8)$&192&2763300&0.360288&&779&3&$(0,0,3)$&72&2762150&0.360138\\
754&3&$(0,0,3)$&72&2763230&0.360279&&780&1&$(0,0,1)$&24&2762140&0.360137\\
755&1&$(0,1,0)$&12&2763220&0.360278&&781&3&$(0,0,3)$&72&2762120&0.360134\\
\hline
\end{tabular}
\end{center}
\end{footnotesize}
\end{sidewaystable}\afterpage{\clearpage}

\setcounter{table}{1}
\renewcommand{\thetable}{A\arabic{table}}
\begin{sidewaystable}[ht]
\caption{The 766 smallest values of the optimal conditional expected return at full pay Jacks or Better, given the initial hand, page~9.  In each case the optimal strategy is to draw five new cards.  See Table A1 for the 387 largest values.\medskip}
\tabcolsep=1.0mm
\catcode`@=\active\def@{\phantom{0}}
\begin{footnotesize}
\begin{center}
\begin{tabular}{rrcrrccrrcrrc}
\hline
\noalign{\smallskip}
no. & no.~of    & numbers of       & prob.       & opt.~condl.          & opt.~condl.  &@@@@@@&  no. & no.~of    & numbers of       & prob.       & opt.~condl.          & opt.~condl. \\
    & equiv.    & equiv.~classes   & ${}\times\binom{52}{5}$ & expected@ & expected     &&      & equiv.    & equiv.~classes   & ${}\times\binom{52}{5}$ & expected@ & expected  \\
    & classes   & of sizes         &             & return@@             & return       &&      & classes   & of sizes         &             & return@@             & return \\
    &           & 4, 12, and 24    &         & ${}\times5\binom{47}{5}@$ & (rounded)    &&      &           & 4, 12, and 24    &         & ${}\times5\binom{47}{5}@$ & (rounded) \\
\noalign{\smallskip}
\hline
\noalign{\smallskip}
804&7&$(0,0,7)$&168&2761430&0.360044&&830&7&$(0,0,7)$&168&2760550&0.359930\\
805&4&$(0,3,1)$&60&2761420&0.360043&&831&2&$(0,0,2)$&48&2760540&0.359928\\
806&3&$(0,0,3)$&72&2761390&0.360039&&832&4&$(0,4,0)$&48&2760500&0.359923\\
807&6&$(0,0,6)$&144&2761350&0.360034&&833&4&$(0,0,4)$&96&2760440&0.359915\\
\hline
\end{tabular}
\end{center}
\end{footnotesize}
\end{sidewaystable}\afterpage{\clearpage}

\setcounter{table}{1}
\renewcommand{\thetable}{A\arabic{table}}
\begin{sidewaystable}[ht]
\caption{The 766 smallest values of the optimal conditional expected return at full pay Jacks or Better, given the initial hand, page~10.  In each case the optimal strategy is to draw five new cards.  See Table A1 for the 387 largest~values.\medskip}
\tabcolsep=1.0mm
\catcode`@=\active\def@{\phantom{0}}
\begin{footnotesize}
\begin{center}
\begin{tabular}{rrcrrccrrcrrc}
\hline
\noalign{\smallskip}
no. & no.~of    & numbers of       & prob.       & opt.~condl.          & opt.~condl.  &@@@@@@&  no. & no.~of    & numbers of       & prob.       & opt.~condl.          & opt.~condl. \\
    & equiv.    & equiv.~classes   & ${}\times\binom{52}{5}$ & expected@ & expected     &&      & equiv.    & equiv.~classes   & ${}\times\binom{52}{5}$ & expected@ & expected  \\
    & classes   & of sizes         &             & return@@             & return       &&      & classes   & of sizes         &             & return@@             & return \\
    &           & 4, 12, and 24    &         & ${}\times5\binom{47}{5}@$ & (rounded)    &&      &           & 4, 12, and 24    &         & ${}\times5\binom{47}{5}@$ & (rounded) \\
\noalign{\smallskip}
\hline
\noalign{\smallskip}
856&15&$(0,0,15)$&360&2759710&0.359820&&882&1&$(0,0,1)$&24&2758480&0.359660\\
857&6&$(0,0,6)$&144&2759700&0.359819&&883&8&$(0,0,8)$&192&2758470&0.359658\\
858&3&$(0,0,3)$&72&2759680&0.359816&&884&4&$(0,0,4)$&96&2758460&0.359657\\
859&5&$(0,4,1)$&72&2759540&0.359798&&885&11&$(0,0,11)$&264&2758430&0.359653\\
\hline
\end{tabular}
\end{center}
\end{footnotesize}
\end{sidewaystable}\afterpage{\clearpage}

\setcounter{table}{1}
\renewcommand{\thetable}{A\arabic{table}}
\begin{sidewaystable}[ht]
\caption{The 766 smallest values of the optimal conditional expected return at full pay Jacks or Better, given the initial hand, page~11.  In each case the optimal strategy is to draw five new cards.  See Table A1 for the 387 largest~values.\medskip}
\tabcolsep=1.0mm
\catcode`@=\active\def@{\phantom{0}}
\begin{footnotesize}
\begin{center}
\begin{tabular}{rrcrrccrrcrrc}
\hline
\noalign{\smallskip}
no. & no.~of    & numbers of       & prob.       & opt.~condl.          & opt.~condl.  &@@@@@@&  no. & no.~of    & numbers of       & prob.       & opt.~condl.          & opt.~condl. \\
    & equiv.    & equiv.~classes   & ${}\times\binom{52}{5}$ & expected@ & expected     &&      & equiv.    & equiv.~classes   & ${}\times\binom{52}{5}$ & expected@ & expected  \\
    & classes   & of sizes         &             & return@@             & return       &&      & classes   & of sizes         &             & return@@             & return \\
    &           & 4, 12, and 24    &         & ${}\times5\binom{47}{5}@$ & (rounded)    &&      &           & 4, 12, and 24    &         & ${}\times5\binom{47}{5}@$ & (rounded) \\
\noalign{\smallskip}
\hline
\noalign{\smallskip}
908&2&$(0,0,2)$&48&2757660&0.359553&&934&7&$(0,0,7)$&168&2756400&0.359388\\
909&1&$(0,0,1)$&24&2757550&0.359538&&935&4&$(0,0,4)$&96&2756390&0.359387\\
910&10&$(0,0,10)$&240&2757510&0.359533&&936&2&$(0,2,0)$&24&2756380&0.359386\\
911&10&$(0,0,10)$&240&2757470&0.359528&&937&9&$(0,0,9)$&216&2756350&0.359382\\
\hline
\end{tabular}
\end{center}
\end{footnotesize}
\end{sidewaystable}\afterpage{\clearpage}

\setcounter{table}{1}
\renewcommand{\thetable}{A\arabic{table}}
\begin{sidewaystable}[ht]
\caption{The 766 smallest values of the optimal conditional expected return at full pay Jacks or Better, given the initial hand, page~12.  In each case the optimal strategy is to draw five new cards.  See Table A1 for the 387 largest~values.\medskip}
\tabcolsep=1.0mm
\catcode`@=\active\def@{\phantom{0}}
\begin{footnotesize}
\begin{center}
\begin{tabular}{rrcrrccrrcrrc}
\hline
\noalign{\smallskip}
no. & no.~of    & numbers of       & prob.       & opt.~condl.          & opt.~condl.  &@@@@@@&  no. & no.~of    & numbers of       & prob.       & opt.~condl.          & opt.~condl. \\
    & equiv.    & equiv.~classes   & ${}\times\binom{52}{5}$ & expected@ & expected     &&      & equiv.    & equiv.~classes   & ${}\times\binom{52}{5}$ & expected@ & expected  \\
    & classes   & of sizes         &             & return@@             & return       &&      & classes   & of sizes         &             & return@@             & return \\
    &           & 4, 12, and 24    &         & ${}\times5\binom{47}{5}@$ & (rounded)    &&      &           & 4, 12, and 24    &         & ${}\times5\binom{47}{5}@$ & (rounded) \\
\noalign{\smallskip}
\hline
\noalign{\smallskip}
960&7&$(0,0,7)$&168&2755640&0.359289&&986&3&$(0,0,3)$&72&2754600&0.359154\\
961&12&$(0,0,12)$&288&2755590&0.359283&&987&2&$(0,0,2)$&48&2754590&0.359152\\
962&1&$(0,0,1)$&24&2755550&0.359278&&988&5&$(0,0,5)$&120&2754470&0.359137\\
963&4&$(0,0,4)$&96&2755500&0.359271&&989&17&$(0,0,17)$&408&2754430&0.359132\\
\hline
\end{tabular}
\end{center}
\end{footnotesize}
\end{sidewaystable}\afterpage{\clearpage}

\setcounter{table}{1}
\renewcommand{\thetable}{A\arabic{table}}
\begin{sidewaystable}[ht]
\caption{The 766 smallest values of the optimal conditional expected return at full pay Jacks or Better, given the initial hand, page~13.  In each case the optimal strategy is to draw five new cards.  See Table A1 for the 387 largest~values.\medskip}
\tabcolsep=1.0mm
\catcode`@=\active\def@{\phantom{0}}
\begin{footnotesize}
\begin{center}
\begin{tabular}{rrcrrccrrcrrc}
\hline
\noalign{\smallskip}
no. & no.~of    & numbers of       & prob.       & opt.~condl.          & opt.~condl.  &@@@@@@&  no. & no.~of    & numbers of       & prob.       & opt.~condl.          & opt.~condl. \\
    & equiv.    & equiv.~classes   & ${}\times\binom{52}{5}$ & expected@ & expected     &&      & equiv.    & equiv.~classes   & ${}\times\binom{52}{5}$ & expected@ & expected  \\
    & classes   & of sizes         &             & return@@             & return       &&      & classes   & of sizes         &             & return@@             & return \\
    &           & 4, 12, and 24    &         & ${}\times5\binom{47}{5}@$ & (rounded)    &&      &           & 4, 12, and 24    &         & ${}\times5\binom{47}{5}@$ & (rounded) \\
\noalign{\smallskip}
\hline
\noalign{\smallskip}
1012&2&$(0,0,2)$&48&2753030&0.358949&&1038&8&$(0,0,8)$&192&2751550&0.358756\\
1013&8&$(0,0,8)$&192&2752990&0.358944&&1039&2&$(0,0,2)$&48&2751470&0.358746\\
1014&7&$(0,0,7)$&168&2752920&0.358935&&1040&7&$(0,0,7)$&168&2751400&0.358737\\
1015&4&$(0,0,4)$&96&2752910&0.358933&&1041&1&$(0,0,1)$&24&2751390&0.358735\\
\hline
\end{tabular}
\end{center}
\end{footnotesize}
\end{sidewaystable}\afterpage{\clearpage}

\setcounter{table}{1}
\renewcommand{\thetable}{A\arabic{table}}
\begin{sidewaystable}[ht]
\caption{The 766 smallest values of the optimal conditional expected return at full pay Jacks or Better, given the initial hand, page~14.  In each case the optimal strategy is to draw five new cards.  See Table A1 for the 387 largest~values.\medskip}
\tabcolsep=1.0mm
\catcode`@=\active\def@{\phantom{0}}
\begin{footnotesize}
\begin{center}
\begin{tabular}{rrcrrccrrcrrc}
\hline
\noalign{\smallskip}
no. & no.~of    & numbers of       & prob.       & opt.~condl.          & opt.~condl.  &@@@@@@&  no. & no.~of    & numbers of       & prob.       & opt.~condl.          & opt.~condl. \\
    & equiv.    & equiv.~classes   & ${}\times\binom{52}{5}$ & expected@ & expected     &&      & equiv.    & equiv.~classes   & ${}\times\binom{52}{5}$ & expected@ & expected  \\
    & classes   & of sizes         &             & return@@             & return       &&      & classes   & of sizes         &             & return@@             & return \\
    &           & 4, 12, and 24    &         & ${}\times5\binom{47}{5}@$ & (rounded)    &&      &           & 4, 12, and 24    &         & ${}\times5\binom{47}{5}@$ & (rounded) \\
\noalign{\smallskip}
\hline
\noalign{\smallskip}
1064&11&$(0,0,11)$&264&2750560&0.358627&&1090&21&$(0,0,21)$&504&2748880&0.358408\\
1065&4&$(0,0,4)$&96&2750550&0.358626&&1091&8&$(0,0,8)$&192&2748790&0.358396\\
1066&1&$(0,0,1)$&24&2750470&0.358615&&1092&9&$(0,0,9)$&216&2748720&0.358387\\
1067&1&$(0,0,1)$&24&2750440&0.358611&&1093&4&$(0,0,4)$&96&2748590&0.358370\\
\hline
\end{tabular}%
\end{center}
\end{footnotesize}
\end{sidewaystable}\afterpage{\clearpage}


\begin{thebibliography}{00}

\bibitem{A07} Alspach, Brian (2007) Enumerating poker hands using group theory.  In Stewart N.~Ethier and William R.~Eadington, editors, \emph{Optimal Play: Mathematical Studies of Games and Gambling}, pages 121--129. Institute for the Study of Gambling and Commercial Gaming, University of Nevada, Reno.

\bibitem{D96} Dancer, Bob (1996) \textit{9-6 Jacks or Better Video Poker: A Complete ``How to Beat the Casino'' Discussion.} Published by the author.  (Available at UNLV Lied Library Special Collections.)

\bibitem{DD03} Dancer, Bob and Daily, Liam W. (2003) \emph{A Winner's Guide to Jacks or Better}. Compton Dancer Consulting Inc., Las Vegas.  (2nd ed., 2004.)

\bibitem{E10} Ethier, Stewart N. (2010) \textit{The Doctrine of Chances: Probabilistic Aspects of Gambling}. Springer, Berlin.\\ \url{http://www.math.utah.edu/~ethier/sample.pdf}.

\bibitem{F89} Frome, Lenny (1989) \emph{Expert Video Poker for Las Vegas}.  Compu-Flyers, Las Vegas.

\bibitem{G96} Gordon, Edward (1996) Analysis of video poker.  \textit{Gaming Research and Review Journal} \textbf{3} (1) 69--80.\\ \url{http://digitalscholarship.unlv.edu/grrj/vol3/iss1/6/}.

\bibitem{G00} Gordon, Edward (2000) An accurate analysis of video poker. In Olaf Vancura, Judy A.~Cornelius, and William R.~Eadington, editors, \emph{Finding the Edge: Mathematical Analysis of Casino Games}, pages 379--392. Institute for the Study of Gambling and Commercial Gaming, University of Nevada, Reno.

\bibitem{K12} Kim, John Jungtae (2012) Optimal strategy hand-rank table for Jacks or Better, Double Bonus, and Joker Wild.  Masters thesis, McMaster University. \url{https://macsphere.mcmaster.ca/bitstream/11375/11863/1/fulltext.pdf}.

\bibitem{M06} Marshall, Marten (2006) \emph{The Poker Code: The Correct Draw for Every Video Poker Hand of Deuces Wild}. Marshall House, Inc.

\bibitem{P92} Paymar, Dan (1992) \emph{Video Poker: Precision Play}. Published by the author.


\bibitem{S15} Shackleford, Michael W. (2016) My methodology for video poker analysis.  \textit{Wizard of Odds} website.\\  \url{http://wizardofodds.com/games/video-poker/methodology/}.

\bibitem{WS92} Weber, Glenn and Scruggs, W. Todd (1992) A mathematical and computer analysis of video poker.  In William R.~Eadington and Judy A.~Cornelius, editors, \emph{Gambling and Commercial Gaming: Essays in Business, Economics, Philosophy and Science}, pages 625--633. Institute for the Study of Gambling and Commercial Gaming, University of Nevada, Reno.

\bibitem{W88} Wong, Stanford (1988) \textit{Professional Video Poker}.  Pi Yee Press, La Jolla, CA.

\end{thebibliography}
\end{document}